%% file: Poisson_Revision_v1.tex
\documentclass[times,11pt]{article} 
\usepackage[small,compact]{titlesec}
\usepackage{verbatim}
\usepackage{wrapfig}
\usepackage{mathrsfs}
\usepackage{amsmath,amssymb,latexsym}
\usepackage{amsfonts}
\usepackage{times,lastpage,bm,xspace}
\usepackage{mathrsfs,bbm}
\usepackage{dsfont}
\usepackage{algorithm}
\usepackage{algorithmic}
\usepackage{stmaryrd}
\usepackage{booktabs}
\pdfoutput=1
\usepackage[pdftex]{graphicx}
\usepackage{url}

\usepackage[usenames]{color}
\definecolor{plum}  {rgb}{.6,0,.6}
\definecolor{forest}  {rgb}{0,.6,0}
\definecolor{midnight}  {rgb}{0,0,.8}
\usepackage[pdftex,plainpages=false,pdfpagelabels,colorlinks=true,urlcolor=midnight,linkcolor=midnight,citecolor=forest]{hyperref}

\usepackage[caption=false,font=footnotesize]{subfig}

\usepackage{amsmath, amsfonts, amssymb, amsthm,color,bm}
\input{AMS_commands}

\begin{document}

\begin{center}
{\bf{\LARGE{Minimax Optimal Convex Methods for Poisson Inverse Problems under $\ell_q$-ball Sparsity}}}

        \vspace*{.1in}
	\begin{tabular}{cc}
	Yuan Li$^{1}$ & Garvesh Raskutti$^{2}$\\
	\end{tabular}

	\vspace*{.1in}

	\begin{tabular}{c}
	  $^1$ Department of Statistics, University of Wisconsin-Madison\\
	  $^2$ Department of Statistics, Computer Science, and Electrical and Computer Engineering, \\ University of Wisconsin-Madison\\
	\end{tabular}

	\vspace*{.1in}

\end{center}

\begin{abstract}
In this paper, we study the minimax rates and provide an implementable convex algorithm for Poisson inverse problems under weak sparsity and physical constraints. In particular we assume the model $y_i \sim \mbox{Poisson}(Ta_i^{\top}f^*)$ for $1 \leq i \leq n$ where $T \in \mathbb{R}_+$ is the intensity, and we impose weak sparsity on $f^* \in \mathbb{R}^p$ by assuming $f^*$ lies in an $\ell_q$-ball when rotated according to an orthonormal basis $D \in \mathbb{R}^{p \times p}$. In addition, since we are modeling real physical systems we also impose positivity and flux-preserving constraints on the matrix $A = [a_1, a_2,...,a_n]^{\top}$ and the function $f^*$. We prove minimax lower bounds for this model which scale as $R_q (\frac{\log p}{T})^{1 - q/2}$ where it is noticeable that the rate depends on the intensity $T$ and not the sample size $n$. We also show that an $\ell_1$-based regularized least-squares estimator achieves this minimax lower bound, provided a suitable restricted eigenvalue condition is satisfied. Finally we prove that provided $n \geq \tilde{K} \log p$ where $\tilde{K} = O(R_q (\frac{\log p}{T})^{- q/2})$ represents an approximate sparsity level, our restricted eigenvalue condition and physical constraints are satisfied for random bounded ensembles. We also provide numerical experiments that validate our mean-squared error bounds. Our results address a number of open issues from prior work on Poisson inverse problems that focuses on strictly sparse models and does not provide guarantees for convex implementable algorithms.
\end{abstract}

\section{Introduction}
\label{SecIntro}

Large-scale Poisson inverse problems arise in a number of applications where counts are modeled using a Poisson distribution. Examples include imaging (see e.g.~\cite{riceCamera,CandesSIM}), conventional fluorescence microscopy (see e.g.~\cite{nimmi_spie,mca,gmca}), network flow analysis (see e.g.~\cite{ElephantsMice,CounterBraids,fishing,expander_pcs}), DNA analysis (see e.g.~\cite{LaureDNA}) and there are many more. In all these problems, a small number of events (e.g. photons hitting a sensor, packets being output, etc.) are observed and these are modeled using a Poisson distribution. In many of these applications, the number of observed events is small relative to the number of model parameters meaning we are in the so-called \emph{high-dimensional} setting.

One standard approach to model the observed counts in the settings above is via a high-dimensional Poisson inverse problem. Specifically, $(y_i)_{i=1}^n$ follows a Poisson distribution and if $y = (y_1, y_2,...,y_n)^{\top} \in \mathbb{R}^n$, we consider the Poisson linear model defined as (see e.g. \cite{dense_pcs, jiang2014minimax, willett2011poisson}): 
\begin{equation}
\label{equation:model}
y\sim\text{Poisson}(TAf^*),
\end{equation}
where $A\in\mathbb{R}^{n\times p}_{+}$ is a sensing matrix corresponding to the $n$ different projections of our signal of interest $f^*\in\mathbb{R}^p_{+}$ and $T\in\mathbb{R}_{+}$ is the known total intensity. In particular, (\ref{equation:model}) is a shorthand expression for the model
$$y_i\sim\text{Poisson}(T\sum_{j=1}^pA_{ij}f^*_j),~~i = 1,...,n,$$
where the $y_i$'s are independent. Here our goal is to learn the underlying parameter vector $f^*$ based on the observed counts $(y_i)_{i=1}^n$ where $A$ is known. Furthermore $p \gg n$ since we are in the high-dimensional setting.

Since we are interested in modeling real physical systems corresponding to the applications described above, additional physical constraints are required on $f^*$ and $A$ as in \cite{dense_pcs, jiang2014minimax, willett2011poisson}. Since $f^*$ corresponds to the rate at which events occur, $f^*\succeq 0$. Further, we impose the normalization $\|f^*\|_1=1$ which is also used in \cite{jiang2014minimax} so that $\|f^*\|_1=1$ which is a measure of signal strength is normalized. In addition $A$ must be composed of non-negative real numbers with each column summing to at most one. Specifically $A$ must satisfy the following physical constraints:
\begin{eqnarray}
\label{eq:constraintA1}
A & \succeq & 0 \\
\label{eq:constraintA2}
A^{\top}\mathbbm{1}_{n \times 1} & \preceq & \mathbbm{1}_{p \times 1}.
\end{eqnarray}
The first constraint \eqref{eq:constraintA1} is referred to as \emph{positivity} that is $A_{ij} \geq 0~\forall (i,j)$ which ensures that $A f^* \succeq 0$ provided $f^* \succeq 0$ and the second constraint \eqref{eq:constraintA2} corresponds to a  \emph{flux-preserving} constraint which ensures that $\|A f^* \|_1 \leq \|f^*\|_1$. Both constraints are natural for the applications described above since counts must be non-negative and the output flux or energy can not exceed the input flux or energy.

Recent work of \cite{jiang2014minimax} provides minimax optimal rates for Poisson inverse problems described here, when $f^*$ is \emph{strictly sparse} in a basis spanned by the columns of the orthonormal matrix $D \in \mathbb{R}^{p \times p}$, meaning that only a small subset of entries of $D^{\top}f^*$ are non-zero. Furthermore, \cite{jiang2014minimax} provides only minimax upper and lower bounds and do not provide theoretical guarantees for an implementable method. 

In many scenarios (e.g. imaging, network flow analysis), the signal of interest $D^{\top}f^*$ is \emph{weakly sparse} meaning that low-dimensional structure is imposed on $D^{\top}f^*$ by requiring that its co-efficients need not be zero, but many co-efficients make a very small contribution to the overall signal. To be precise we write $D=[d_1,...,d_p]$ for $d_j\in\mathbb{R}^p$ for all $j$ and assume that $d_1=p^{-1/2}\mathbbm{1}_{p \times 1}$ and denote $\bar{D}= [d_2,...,d_p]\in\mathbb{R}^{p\times(p-1)}$. We define $\theta^*=D^{\top}f^*\in\mathbb{R}^p$, and by construction $\theta^*_1=1/\sqrt{p}$. To impose weak sparsity, we require the signal $\bar{\theta}^*\triangleq\bar{D}^{\top}f^*\in\mathbb{R}^{p-1}$ to lie in an $\ell_q$-ball, meaning $\|\bar{\theta}^*\|_q^q=\|\bar{D}^{\top}f^*\|_q^q=\sum_{j=1}^{p-1}|(\bar{D}^{\top}f^*)_j|^q \leq R_q$ where $0 < q \leq 1$. In this paper, we study the Poisson model \eqref{equation:model} under the positivity and flux-preserving constraints and $\ell_q$-ball sparsity.

To summarize, we assume $f^*$ belongs to the following set:
$$\mathcal{F}_{p,q,D}=\{f\in\mathbb{R}^p_+~:~\|f\|_1=1,~\|\bar{D}^{\top}f\|^q_q\leq R_q\}.$$
Note that our $\ell_q$-ball assumption ensures that many of the co-efficients are small and the convention that is often used is $q=0$ corresponds to the strictly sparse case studied in \cite{jiang2014minimax}. In this paper, we study minimax rates for the mean-squared $\ell_2$-error for the Poisson inverse problem \eqref{equation:model} where $f^*$ lies in $\mathcal{F}_{p,q,D}$ where $0 < q \leq 1$. That is we provide (1) a lower bound with high probability on the following quantity:
$$\min_{f(A, y)}\max_{f^{*}\in\mathcal{F}_{p,q,D}}\|f-f^{*}\|_2^2,$$
where the minimum is taken over measurable functions of $(A, y)$; and (2) we show that a convex $\ell_1$-penalized approach achieves this optimal rate. As discussed in \cite{foucart2010gelfand}, the geometry of $\ell_q$-balls means only when $0 < q \leq 1$ we can achieve desirable mean-squared error by using sparse vector approximation. Hence we focus on $0<q\leq1$.

\subsection{Our Contributions}
Our paper makes the following novel contributions:
\begin{itemize}
\item Provide a minimax lower bound which scales as $R_q \big(\frac{\log p}{T} \big)^{1 - q/2}$.
\item Show that our minimax lower bound can be achieved by an $\ell_1$-based convex method under a suitable restricted eigenvalue condition.
\item Prove that random bounded ensembles satisfy the restricted eigenvalue condition along with the imposed physical constraints provided $n \geq \tilde{K} \log p$ where $\tilde{K} = O(R_q \big(\frac{\log p}{T} \big)^{-q/2})$ represents the effective sparsity.
\item Provide a simulation study that supports our theoretical mean-squared error.
\end{itemize}

Our bounds are consistent with the intuition from \cite{jiang2014minimax} under strictly sparse models, where the intensity $T$ and not the sample size $n$ influences the minimax rate. To further support this intuition we provide a comparison of our result to the linear Gaussian inverse problem studied in \cite{RasWaiYu11}, and show how the minimax rates match the linear Gaussian rate when we set the noise variance $\sigma^2$ in terms of $n$ and $T$ in Section \ref{SecDiscTvsN}.

We point out that it is not straightforward to adapt the techniques developed in \cite{jiang2014minimax} to the $\ell_q$-ball case. The techniques we use involve combining techniques for proving minimax rates in the high-dimensional Gaussian linear inverse problems under weak sparsity used in \cite{RasWaiYu11} and theoretical results for convex implementable methods developed in \cite{Neg10} to the linear Poisson setting. A number of technical challenges arise in analyzing the Poisson inverse problem setting since the noise is now signal-dependent. In particular, to use techniques from \cite{Neg10} in the Poisson setting we need to use and develop two-sided concentration bounds for Poisson inverse problems which build on prior work in \cite{bobkov1998modified}. We go into greater detail on the technical challenges in Sections \ref{SecMain}, \ref{SecDiscussion} and \ref{SecProofs}.

The remainder of this paper is organized as follows: In Section \ref{SecMain} we provide our main assumptions and theoretical results which includes a minimax lower bound, an upper bound for convex methods and discuss matrices $A$ that satisfy the assumptions leading to minimax rates; in Section \ref{SecDiscussion} we discuss implications of our results in particular, comparisons to the linear Gaussian model studied in \cite{RasWaiYu11} and a comparison to the strictly sparse case in \cite{jiang2014minimax}. Numerical experiments are provided in Section \ref{sc:simulation} and proofs are provided in Section \ref{SecProofs}.

\section{Assumptions and Main Results}

\label{SecMain}
In this section, we present our assumptions and main results, which includes a minimax lower bound, an upper bound for a convex $\ell_1$-based approach that matches the minimax lower and finally we show that if $A$ is a random matrix with suitably bounded entries, it satisfies the statistical conditions and physical constraints.

\subsection{Minimax Lower Bound}
\label{SecLower}
We begin by introducing the assumptions for the minimax lower bound.

\bas
\label{as:construction}
There exists constants $a_{\ell}$ and $a_u$ such that $a_{\ell}<a_u$ and a matrix $\tilde{A}\in [a_{\ell}/\sqrt{n},a_u/\sqrt{n}]^{n\times p}$
\begin{equation}
\label{eq:construction}
A=\frac{\tilde{A}+\frac{a_u-2a_{\ell}}{\sqrt{n}}\mathbbm{1}_{n\times p}}{2(a_u-a_{\ell})\sqrt{n}}.
\end{equation}
\eas
Assumption \ref{as:construction} is originally imposed in \cite{jiang2014minimax} and ensures that the positivity and flux-preserving conditions are satisfied. Further, Assumption~\ref{as:construction} ensures the values within the matrix are approximately uniform. By Lemma 2.1 in \cite{jiang2014minimax} if matrix $\tilde{A}$ satisfies $\tilde{A}_{i,j}\in\frac{1}{\sqrt{n}}[a_{\ell},a_u]$ for all entries, then the sensing matrix $A$ will satisfy the physical constraints (\ref{eq:constraintA1}) and (\ref{eq:constraintA2}).

\bas
\label{as:rip}
For all $u\in\mathbb{R}^p,~\|u\|_0\leq2\tilde{K}$, there exists a constant $\delta_{\tilde{K}}\equiv \delta_{\tilde{K}}(n,p)>0$ such that:
$$\|\tilde{A}Du\|^2_2\leq(1+\delta_{\tilde{K}})\|u\|^2_2,$$
where $\tilde{K}$ is a constant satisfies $\tilde{K}=O(R_q(\frac{\log p}{T})^{-\frac{q}{2}}),$ without loss of generality we can assume $\tilde{K}$ is an integer.
\eas
\brem
\label{rm:rip}
Assumption \ref{as:rip} is an upper restricted isometry property condition similar to that imposed in \cite{jiang2014minimax}. The main difference is that we use a different sparsity parameter $\tilde{K}$ which depends on $n$, $p$, and $T$ instead of $s$ in the strict sparsity case. This assumption will be used for both the lower bound and upper bound for the $\ell_1$-based method. As pointed out in \cite{jiang2014minimax}, Assumption \ref{as:rip} holds with $0<\delta_{\tilde{K}}<1$ which occurs if $n\geq \tilde{K}\log p$ for a re-scaled Bernoulli ensemble matrix $\tilde{A}$ with $P(\tilde{A}_{ij}=\frac{1}{\sqrt{n}})=P(\tilde{A}_{ij}=-\frac{1}{\sqrt{n}})=\frac{1}{2}$,  with probability at least $1-e^{-C_1n}$ using results in \cite{Baraniuk08}.
\erem

%

Finally we define an $s$-sparse localization quantity also introduced in \cite{jiang2014minimax}. The interaction between the orthonormal basis matrix $D$ and the sparsity constraint has an effect on the lower bound which is captured by this $s$-sparse localization quantity:
\bde
$\lambda_s$ is said to be the $s$-sparse localization quantity of a matrix $X$ if 
$$\lambda_s=\lambda_s(X):=\max_{v\in\{-1,0,1\}^p;\|v\|_0=s}\|Xv\|_{\infty}.$$
\ede
Our minimax lower bound depends on $\lambda_{k}(\bar{D})$. Different scaling for the $k$-sparse localization constant $\lambda_k(\bar{D})$ with basis $D$ for both Fourier and wavelet transforms are provided in \cite{jiang2014minimax}. Now we present the minimax lower bound.

\btheos
\label{theorem:lower}
If $f^{*}\in\mathcal{F}_{p,q,D}$, $p \geq \max(260,\frac{33\tilde{K}}{2}+1)$, Assumption \ref{as:construction} to \ref{as:rip} hold with $0<\delta_{\tilde{K}}<1$, let $\lambda_{k}=\lambda_k(\bar{D})$ be the $k$-$sparse$ localization quantity of $\bar{D}$, then if $\tilde{K}=O(R_q(\frac{\log p}{T})^{-\frac{q}{2}})$,
there exists a constant $C_L>0$ that depends on $a_u, a_{\ell}~\text{and}~\delta_{\tilde{K}}$ such that
\begin{eqnarray}
\label{eq:lowerbound}
\min_{f}\max_{f^{*}\in\mathcal{F}_{p,q,D}}\|f-f^{*}\|_2^2\geq C_L \max_{1 \leq k \leq \tilde{K}}\left\{ \min\left(\frac{k}{p^2\lambda_{k}^2},\frac{k\log\frac{p}{k}}{T}\right)\right \}
\end{eqnarray}
with probability greater than $\frac{1}{2}$.
\etheos

{\bf Remarks:}
\begin{itemize}	
\item Note that in Theorem \ref{theorem:lower} when $k=\tilde{K}=O(R_q(\frac{\log p}{T})^{-\frac{q}{2}})$, the second term $\frac{k\log\frac{p}{k}}{T}$ in the lower bound scales as $R_q(\frac{\log p}{T})^{1-\frac{q}{2}}$.

\item In the case $q = 0$, the minimax rate is $\frac{s \log p}{T}$ as proven in \cite{jiang2014minimax}. Note that if we set $s = s_q = R_q (\frac{\log p}{T})^{-q/2}$ leading to the overall rate $R_q (\frac{\log p}{T})^{1-q/2}$. This interpretation is consistent with the case of Gaussian linear models discussed in \cite{RasWaiYu11}.

\item Note that as in the case $q = 0$ discussed in \cite{jiang2014minimax}, the minimax lower bound depends on the intensity $T$ and not the sample size $n$. This may initially seem counter-intuitive since the sample size $n$ plays no role in the minimax lower bound. This phenomenon arises due to the signal-dependent noise and the flux-preserving constraint which we discuss in greater detail in Section \ref{SecDiscTvsN}. But note that Assumption \ref{as:rip} only holds with $0<\delta_{\tilde{K}}<1$ which typically depends on $n$. For example if $\tilde{A}$ is a re-scaled Bernoulli ensemble matrix, we need $n\geq C\tilde{K}\log p$ for a constant $C>0$ to ensure that $0<\delta_{\tilde{K}}<1$.

\item Values for $\lambda_k$ are displayed in Table 1 of \cite{jiang2014minimax} for the discrete cosine transform (DCT), discrete Hadamard transform (DHT), and a discrete Haar wavelet basis (DWT). In particular for the DCT and DHT basis, $\lambda_k = \frac{\sqrt{2}k}{\sqrt{p}}$ and for DWT, $\lambda_k = \frac{1}{\sqrt{2}-1}$. Therefore the first term $\frac{k}{p^2\lambda_k^2}$ in the lower bound will be $\frac{1}{pk}$ for DCT and DHT, and $\frac{k}{p^2}$ for DWT. Since when $k=\tilde{K}$ the second term $\frac{k\log\frac{p}{k}}{T}$ scales as $R_q(\frac{\log p}{T})^{1-\frac{q}{2}}$, then for DCT and DHT basis, $R_q(\frac{\log p}{T})^{1-\frac{q}{2}}$ will be smaller than the first term $\frac{1}{p\tilde{K}}$ when $p^{\frac{1}{1-q}}\preceq T$, and for DWT, $R_q(\frac{\log p}{T})^{1-\frac{q}{2}}$ will be smaller than $\frac{\tilde{K}}{p^2}$ when $p^2\preceq T$.

\item Since we require that $p\geq\frac{33\tilde{K}}{2}+1$, this leads to $p^{\frac{2}{q}}\succeq T$. Thus by combining with above discussion we can see that for DCT and DHT, the lower bound scales as $R_q(\frac{\log p}{T})^{1-\frac{q}{2}}$ when we are in the setting $p^{\frac{1}{1-q}}\preceq T\preceq p^{\frac{2}{q}}$; and for DWT the lower bound scales as $R_q(\frac{\log p}{T})^{1-\frac{q}{2}}$ when we are in the setting $p^2\preceq T\preceq p^{\frac{2}{q}}$.

\item Although many steps of the proof are similar to the strictly sparse case in \cite{jiang2014minimax}, the $\ell_q$-ball sparsity is more challenging than the strictly sparse case since the $\ell_q$-ball is a compact set. In the proof for Theorem~\ref{theorem:lower} we construct a packing set for the intersection of the $\ell_q$-ball with the physical constraints on $f^*$. Our packing set is based on a combination of the hypercube construction provided in \cite{Kuh01} along with the construction in \cite{jiang2014minimax} which incorporates the positivity and flux-preserving constraints. Further details are provided in Section \ref{ProofLower}.

\end{itemize}

\subsection{$\ell_1$-based Method}
\label{Secl1method}
In this section we present an $\ell_2$-error upper bound for an $\ell_1$-based estimator by adapting existing results and techniques in \cite{Neg10} to our Poisson inverse problem setting. The estimator we consider is the standard Lasso estimator:
$$\hat{\theta}_{\lambda_n}\in\text{arg}\min_{\theta\in\mathbb{R}^p,\theta_1=\frac{1}{\sqrt{p}}}\frac{1}{n}\|\frac{n}{T}(y-TAD\theta)\|_2^2+\lambda_n\|\theta\|_1$$
or equivalently
\begin{equation}
\label{equation:l1problem}
\hat{\theta}_{\lambda_n}\in\text{arg}\min_{\theta\in\mathbb{R}^p,\theta_1=\frac{1}{\sqrt{p}}}\frac{n}{T^2}\|y-TAD\theta\|_2^2+\lambda_n\|\theta\|_1,
\end{equation}
where $\lambda_n>0$ is the regularization parameter. We can then get the estimator $\hat{f}_{\lambda_n}$ for $f^*$ by $\hat{f}_{\lambda_n}=D\hat{\theta}_{\lambda_n}$. Next we introduce two further assumptions:

\bas
\label{Assumption:ntildek}
There exists a constant $c_0>0$ such that $n\geq c_0\tilde{K}\log p$.
\eas
\brem
Assumption \ref{Assumption:ntildek} will be used to derive the so-called restricted eigenvalue condition (see e.g. \cite{BicRitTsy08, RasWaiYu10b, GeerBuhl09}) for matrix $A\bar{D}$. In Section \ref{SecLower} we can see that this assumption also ensures that Assumption \ref{as:rip} hold with $0<\delta_{\tilde{K}}<1$ for re-scaled Bernoulli ensemble matrix. It is important to note that this assumption is equivalent to assume that 
\begin{eqnarray}
\label{eq:upperT}
T<M_1(\frac{n}{R_q\log p})^{\frac{2}{q}}\log p
\end{eqnarray}
for some positive constant $M_1$.
\erem
\bas
\label{Assumption:L1RE}
There are strictly positive constants $(k_1,k_2)$ that depend on $a_u,~a_{\ell}$ such that
$$\sqrt{n}\|A\bar{D}x\|_2\geq k_1\|x\|_2-k_2\sqrt{\frac{\log p}{n}}\|x\|_1,~~\forall x\in\mathbb{R}^{p-1}.$$
\eas

Assumption \ref{Assumption:L1RE} is also used to derive the restricted eigenvalue condition for the matrix $A\bar{D}$. This assumption holds for many appropriate choices of A, as we show in Theorem \ref{Theorem:RE} in Section \ref{SecRECondition}.\\ 

The upper bound is as follows:
\btheos
\label{ThmL1}
If $f^{*}\in\mathcal{F}_{p,q,D}$, $T>2n\log p$ and Assumption \ref{as:construction} to \ref{Assumption:L1RE} hold with $0<\delta_{\tilde{K}}<1$. If we choose $\lambda_n=2\sqrt{\frac{32M\log p}{T}}$ with $M=\frac{1+\delta_{\tilde{K}}}{4(a_u-a_{\ell})^2}$, then there exists a constant $C_U>0$ that depends on $a_u, a_{\ell}~\text{and}~\delta_{\tilde{K}}$ such that
\begin{eqnarray}
\label{eq:upperbound}
\|\hat{f}_{\lambda_n}-f^{*}\|_2^2\leq C_UR_q(\frac{\log p}{T})^{1-\frac{q}{2}}
\end{eqnarray}
with probability at least $1-\frac{2}{p-1}$.
\etheos

{\bf Remarks:}
\begin{itemize}
\item Theorem \ref{ThmL1} shows that the upper bound result for this $\ell_1$-based estimator achieves rates $R_q(\frac{\log p}{T})^{1-\frac{q}{2}}$ for $T$ satisfies (\ref{eq:upperT}) and $T>2n\log p$. Since $T$ controls the signal to noise ratio, $T>2n\log p$ is equivalent to ensure the signal to noise ratio is big enough for good estimation.
\item By combining discussion in Section \ref{SecLower} about conditions for $T$ we can see that:
\begin{itemize}
\item If the orthonormal matrix $D$ is DCT or DHC basis, the upper bound and lower bound results match with a mean-squared error rate $R_q(\frac{\log p}{T})^{1-\frac{q}{2}}$ when $p^{\frac{1}{1-q}}\preceq T\preceq (\frac{n}{R_q\log p})^{\frac{2}{q}}\log p$.
\item If $D$ is DWT basis, the upper bound and lower bound results match with a mean-squared error rate $R_q(\frac{\log p}{T})^{1-\frac{q}{2}}$ when $p^2\preceq T\preceq (\frac{n}{R_q\log p})^{\frac{2}{q}}\log p$.
\end{itemize}

\item Though Theorem 2 shows that this upper bound depends explicitly on the intensity $T$ but not the number of observations $n$, it is crucial to know that $n$ actually plays an important role through both Assumption 2.3 and 2.4. It also enters through the required normalization on $A$ in order to guarantee flux-preserving constraint.



\end{itemize}

\subsection{Restricted Eigenvalue Condition}

\label{SecRECondition}
In this section we show that Assumption \ref{Assumption:L1RE} is satisfied with high probability for many choices of random matrices $A$ under the appropriate scaling. In particular we show that the restricted eigenvalue condition is satisfied by matrices $A$ with independent sub-Gaussian entries which include independent Bernoulli ensembles that also satisfy our flux-preserving and positivity constraints.

To characterize the sub-Gaussian parameter of a random variable, we define the Orlicz norm $\|.\|_{\psi_2}$ for a random variable $X \in \mathbb{R}$ as follows:
$$\|X\|_{\psi_2}:=\inf\{t:~\mathbb{E}\exp(X^2/t^2)\leq 2\}.$$
The Orlicz norm as defined above is known to represent the sub-Gaussian parameter of a random variable. For example if $X \sim \mathcal{N}(0,\sigma^2)$, $\|X\|_{\psi_2} = \sigma$. Now we provide a definition of isotropic random vectors introduced in \cite{mendelson2007reconstruction} and \cite{zhou2009restricted}.
\bde[\cite{zhou2009restricted}, Definition 1.3]
Let $Y$ be a random vector in $\mathbb{R}^p$; $Y$ is called isotropic if for every $y\in\mathbb{R}^p$, $\mathbb{E}|\langle Y,y\rangle|^2=\|y\|_2^2$, and is $\psi_2$ with a constant $\alpha$ if for every $y\in\mathbb{R}^p$:
$$\|\langle Y,y\rangle\|_{\psi_2} \leq\alpha \|y\|_2.$$
\ede
Important examples of isotropic vectors are the Gaussian random vector $Y=(h_1,...,h_p)$ where $h_i,~\forall i$ are independent $N(0,1)$ random variables where $\alpha = 1$, and random vectors $Y=(\epsilon_1,...,\epsilon_p)$ where $\epsilon_i,~\forall i$ are independent, symmetric $\pm 1$ Bernoulli random variables also with $\alpha =1$. Now we are ready to state the main theorem for this section:
\btheos
\label{Theorem:RE}
There exists positive constants $c',c''$ for which the following holds. Let $\mu$ to be an isotropic $\psi_2$ probability measure with constant $\alpha\geq1$. Let $X_1,...,X_n\in\mathbb{R}^p$ be independent, distributed according to $\mu$ and define $\Gamma=\sum_{i=1}^n\langle X_i,~.~\rangle e_i$, where $e_i\in\mathbb{R}^n$ is a vector with $i^{\text{th}}$-location to be 1 and all the other locations to be 0. Then with probability at least $1-c'\exp(-c''n)$, for all $x\in\mathbb{R}^p$ we will have
$$\frac{\|x\|_2}{4}-C_\alpha\sqrt{\frac{\log p}{n}}\|x\|_1\leq\frac{\|\Gamma x\|_2}{\sqrt{n}},$$
where $C_\alpha$ is a positive constant only depends on $\alpha$.
\etheos
{\bf Remarks:}
\begin{itemize}
\item Theorem \ref{Theorem:RE} shows that the restricted eigenvalue condition holds for matrices with random sub-Gaussian entries which include both Gaussian and bounded random variables. The proof techniques are based on a combination of techniques from \cite{RasWaiYu10b} for random Gaussian matrices with techniques from \cite{mendelson2007reconstruction} for sub-Gaussian random variables. The proof is provided in Section \ref{SecProofRE}.
\item Based on this theorem, there are many choices of $A$ which satisfy Assumption \ref{Assumption:L1RE}. In our particular context, we also require $A$ to satisfy Assumption \ref{as:construction} so that it satisfies our physical constraints. Hence we require the entries of $A$ to be bounded, and we provide a concrete example below.
\item Theorem \ref{Theorem:RE} is more general than the restricted eigenvalue condition for strictly sparse vectors proven by \cite{zhou2009restricted}. Our result easily adapts to weak $\ell_q$-ball sparse vectors and in addition our result applies to any $x$ that may be random which we address using a peeling argument in our proof.
\end{itemize}

To construct a random matrix $A$ that satisfies the restricted eigenvalue condition and Assumption \ref{as:construction}, let $\tilde{A}$ have the following entries:
$$\mathbb{P}(\tilde{A}_{ij})=\begin{cases}
\frac{1}{2} & \tilde{A}_{ij}=-\sqrt{\frac{1}{n}}\\
\frac{1}{2} & \tilde{A}_{ij}=\sqrt{\frac{1}{n}},\\
\end{cases}$$
then $\sqrt{n}\tilde{A}$ will satisfy the conditions for $\Gamma$ in Theorem \ref{Theorem:RE}. Since we want our result to apply after we apply an orthonormal basis $D$ we use the following Lemma:
\blems
\label{lemma:ADsubgaussian}
Let $\mu$ to be an isotropic $\psi_2$ probability measure with constant $\alpha\geq 1$. And let $X\in\mathbb{R}^p$ be distributed according to $\mu$, then $X\bar{D}\in\mathbb{R}^{p-1}$ is distributed according to another isotropic $\psi_2$ probability measure $\mu'$ with some constant $\alpha'\geq 1$.
\elems

Thus by Lemma \ref{lemma:ADsubgaussian}, $\sqrt{n}\tilde{A}\bar{D}$ satisfies the restricted eigenvalue condition from Theorem \ref{Theorem:RE} and with high probability:
$$\|\tilde{A}\bar{D}x\|_2=\frac{\|\sqrt{n}\tilde{A}\bar{D}x\|_2}{\sqrt{n}}\geq\frac{\|x\|_2}{4}-C\sqrt{\frac{\log p}{n}}\|x\|_1,~~\forall{x}\in\mathbb{R}^{p-1},$$
where $C>0$ is some absolute constant. Note that by the construction of $A$ in Assumption \ref{as:construction} and definition of $\bar{D}$ we have
$$A\bar{D}x=\frac{\tilde{A}+\frac{a_u-2a_{\ell}}{\sqrt{n}}\mathbbm{1}_{n\times p}}{2(a_u-a_{\ell})\sqrt{n}}\bar{D}x=\frac{\tilde{A}\bar{D}x}{2(a_u-a_{\ell})\sqrt{n}}.$$
Then
$$\sqrt{n}\|A\bar{D}x\|_2\geq\frac{\|x\|_2}{8(a_u-a_{\ell})}-\frac{C}{2(a_u-a_{\ell})}\sqrt{\frac{\log p}{n}}\|x\|_1,~~\forall{x}\in\mathbb{R}^{p-1},$$
which satisfies Assumption \ref{Assumption:L1RE}.

\section{Discussion}
\label{SecDiscussion}

In this section, we discuss some of the consequences and intuition for our three main results. In particular we discuss the dependence of the rates on $n$ and $T$ and discuss connections to the results for the Gaussian linear model in \cite{RasWaiYu11}, how the results for the $\ell_q$ case relate to the strictly sparse case developed in \cite{jiang2014minimax} and finally we compare our upper bound to the upper bounds developed in the recent work of \cite{Jiang2015data} based on the weighted Lasso. 

\subsection{Dependence on $T$ and $n$}

\label{SecDiscTvsN}

One of the interesting and perhaps surprising aspects about both the upper and lower bounds is that they depend explicitly on the the intensity $T$ and not on the sample size $n$, aside from the conditions on design. This phenomenon also occurred in the strictly sparse case in \cite{jiang2014minimax} where the rate is $\frac{s \log p}{T}$. To understand this, we relate our rate of $R_q(\frac{\log p}{T})^{1-q/2}$ to the earlier results developed in \cite{RasWaiYu11} for the Gaussian linear model and see how the signal-dependent noise and physical constraints ensure the minimax rate depends on $T$ and not $n$.

In the Gaussian linear model under the $\ell_q$-ball constraint studied in \cite{RasWaiYu11}, we have
$$y=\bar{A}f^*+w,$$
where $y\in\mathbb{R}^n$, $\bar{A}\in\mathbb{R}^{n\times p}$ with $p>n$ and $w\sim \mbox{N}(0,\sigma^2 I_{n \times n})$, and we have the constraint $\|f^*\|_q^q\leq R_q$ with $0<q\leq 1$. \cite{RasWaiYu11} shows that the minimax rate is:
$$\min_{\hat{f}}\max_{\|f^*\|_q^q\leq R_q}\|\hat{f}-f^*\|^2_2\asymp R_q(\frac{\sigma^2 \log p}{n})^{1-\frac{q}{2}},$$
with high probability. In particular, take note of the role of $\sigma^2$ in the minimax rate. Later work by \cite{Neg10} proves that the Lasso estimator achieves this minimax rate. We will show how the dependence of the scaling on $T$ and not $n$ follows from our comparison to the rates for the Gaussian linear model and the impact of $\sigma^2$.

For our Poisson inverse problem we can express the model as follows:
$$y\sim\text{Poisson}(TAf^*),$$ which can be expressed equivalently
$$y=TAf^*+\omega,$$
or
$$y_i = T\sum_{j=1}^pA_{ij}f^*_j + \omega_i,~1\leq i\leq n,$$
where
$$\mathbb{E}(y_i)=T\sum_{j=1}^pA_{ij}f^*_j=\frac{\alpha T}{n},~~1\leq i\leq n,$$
and
$$\mathbb{E}(\omega_i)=0,~~\mbox{Var}(\omega_i)=\mbox{Var}(y_i)=\mathbb{E}(y_i)=\frac{\alpha T}{n},~~1\leq i\leq n,$$
where $\frac{1}{2}\leq\alpha\leq 1$ by Lemma \ref{lemma:Afbound} in section \ref{ProofLower}. Note that this scaling of $A$ follows from the flux-preserving constraint $\|A f\|_1 \leq \|f\|_1$. Since we have $\mathbb{E}(y_i)$ scaling as $\frac{T}{n}$, to ensure that the mean of our observations has the same scaling as in the Gaussian linear model independent of $n$ and $T$, we consider the normalized responses:
$$\tilde{y_i}=\frac{n}{T}y_i=n\sum_{j=1}^pA_{ij}f^*_j+\frac{n}{T}\omega_i,~1\leq i\leq n,$$
where now $\mathbb{E}(\tilde{y}_i)=\alpha$ has a scaling independent of $n$ and $T$, by defining $\tilde{\omega}_i=\frac{n}{T}\omega_i$ we also have
$$\mathbb{E}(\tilde{\omega}_i)=0~~\text{and}~~\mbox{Var}(\tilde{\omega}_i)=\frac{n^2}{T^2}\frac{\alpha T}{n}=\frac{\alpha n}{T}.$$
Hence the combination of the signal-dependent noise and the flux-preserving constraint mean that we have $\sigma^2$ scaling as $\frac{n}{T}$ in the appropriately normalized model. 

Recall that for the Gaussian model, the minimax rate is $R_q(\frac{\sigma^2\log p}{n})^{1-q/2}$, and if we replace $\sigma^2$ by $\frac{\alpha n}{T}$ we will get the minimax rate scaling as $R_q(\frac{\alpha \log p}{T})^{1-q/2}$, which is consistent with our minimax lower bound. Hence, if we want to relate our physically constrained Poisson model to the Gaussian linear model, we need to consider a model with variance $\sigma^2 \propto \frac{n}{T}$.

This observation that the minimax rate depends on the signal intensity $T$ rather than the sample size $n$ was also made in the recent work of \cite{jiang2014minimax}. Our analysis shows that this observation carries over to the $\ell_q$-ball setting. The caveat is that $n$ is required to be sufficiently large to ensure that the restricted eigenvalue is satisfied which is also required in the strictly sparse case.

\subsection{Results without flux-preserving constraint}
\label{SecEllOne}

One of the main contributions of the paper is to develop mean-squared error bounds under the physical constraints, namely the flux-preserving and non-negativity constraint. A natural question to consider is how the results would change if we remove the flux-preserving constraint. In this section we briefly discuss how the error bound changes if we don't have the flux-preserving constraint for matrix $A$. That means we consider a model similar to \cite{Jiang2015data}:
\begin{eqnarray}
y\sim\text{Poisson}(Af^*),
\end{eqnarray}
where $y\in\mathbb{R}_{+}^n$, $A\in\mathbb{R}_{+}^{n\times p}$ and $f\in\mathbb{R}^{p}_{+}$ and no flux-preserving constraint for $A$, instead we assume that $A_{ij}=O(1)$ for $1\leq i\leq n$ and $1\leq j\leq p$. Hence we also remove the normalization condition that $\|f^*\|_1=1$. Without the flux-preserving constraint it is straightforward to show a more general upper bound $\|\hat{f}-f^*\|_2^2\preceq R_q(\frac{\|f^*\|_1\log p}{n})^{1-q/2}$ for arbitrary $\|f^*\|_1$. Hence, if we maintain the condition $\|f^*\|_1=1$ we can show $\|\hat{f}-f^*\|_2^2\preceq R_q(\frac{\log p}{n})^{1-q/2}$. This is consistent with the discussion in Section \ref{SecDiscTvsN} since now the noise variance $\sigma^2$ for observation $y_i$ ($1\leq i\leq n$) is just $\sum_{j=1}^pA_{ij}f_j=O(\|f^*\|_1)$. Similar upper bound results for strictly sparse $f^*$ are discussed in Section 4.5 of \cite{Jiang2015data}. In particular they show $\|\hat{f}_{Lasso}-f^*\|_2^2\preceq\frac{\|f^*\|_1s\log p}{n}$ if $f^*$ contains only $s$ nonzero elements, note that this upper bound also depends on $\|f^*\|_1$.

\subsection{Comparison to related results}

For the strictly sparse case studied by \cite{jiang2014minimax}, the minimax rate scales as $\frac{s \log p}{T}$ whereas in the weakly sparse case in this paper, the minimax rate scales as $R_q(\frac{\log p}{T})^{1-q/2}$. Another way to interpret our result for the $\ell_q$-ball case is that the minimax rate scales as $\frac{s_q \log p}{T}$ where $s_q = O(R_q(\frac{\log p}{T})^{-q/2})$. This can be explained in terms of a bias-variance trade-off to determine how many co-ordinates of $f^*$ should be included in the model and using the $\ell_q$-ball constraint, selecting $s_q = O(R_q(\frac{\log p}{T})^{-q/2})$ with largest magnitude optimizes the bias-variance trade-off to minimize the mean-squared error. This interpretation is used at several points in the proofs of both the minimax lower bound and upper bound. This phenomenon was also observed in the Gaussian linear model case in \cite{RasWaiYu11}.

Another recent related work is by \cite{Jiang2015data} which provides analysis for a weighted Lasso estimator. In \cite{Jiang2015data} sparse Poisson inverse problems under the model $Y\sim \mbox{Poisson}(Af^*)$ are discussed and \cite{Jiang2015data} provides a weighted Lasso estimator $\hat{f}^{WL}$ based on the minimizer of the following optimization problem:
\begin{equation*}
\hat{f}^{WL} \in \arg \min_f \|\tilde{Y} - \tilde{A}f\|_2^2 + \gamma\sum_{j=1}^{p}{d_j |f_j|},
\end{equation*}
where $\tilde{Y}$ and $\tilde{A}$ are shifted and scaled versions of $Y$ and $A$, $\gamma>2$ is a constant and positive weights $(d_j)_{j=1}^p$ are chosen in a specific way to minimize mean-squared error. In the case where all the weights are the same $d_j = \lambda_n$ for all $j$ which corresponds to the ordinary Lasso estimator we analyze in this paper. We summarize their result and show that it is sub-optimal for $\ell_q$-balls. To be clear, the focus of the results in \cite{Jiang2015data} is the strictly sparse case in a number of more general settings than this paper where they achieve optimal or near-optimal mean-squared error. However, they do have a result for approximately sparse models which is not focussed specifically on the $\ell_q$-ball sparsity setting. To summarize their result, they introduce a bias term:
$$B_s :=\max\{\|\tilde{A}(f^*-f_s^*)\|_2^2,\|f^*-f^*_s\|_1\},$$ 
where $s>0$ is an integer and $f_s^*$ is the best $s$-sparse approximation to $f^*$, then they state that
\begin{equation}
\label{equation:wlasso}
\|f^* -\hat{f}^{WL}\|_2^2 \asymp B_s + \lambda^2 s + (1 + 1/\lambda)^2 B_s^2.
\end{equation}
It must be pointed out that \cite{Jiang2015data} analyze a broader choice of weights $(d_j)_{j=1}^p$ but we were unable to find different choices of weights that provided a sharper upper bound in our $\ell_q$-ball context. As discussed in \cite{Jiang2015data}, their analysis yields the optimal rates if $f^*$ is $s$-sparse and the bias term $B_s = 0$. In the case of $\ell_q$-ball sparsity as discussed earlier an appropriate choice for $s$ is $s = s_q = O(R_q(\frac{\log p}{T})^{-\frac{q}{2}})$ and $\lambda_n = O(\sqrt{\frac{\log p}{T}})$, for the $B_s$ term, $\|f^*-f^*_s\|_1$ will be bounded by $R_q\lambda_n^{1-q}$ by inequality (\ref{equation:theta}). By replacing these terms in (\ref{equation:wlasso}), $\|f^*-\hat{f}^{WL}\|^2_2$ will be bounded by $O(R_q\lambda_n^{1-q}+R_q^2\lambda_n^{-2q})$, which scales as $R_q(\frac{\log p}{T})^{\frac{1-q}{2}}+R_q^2(\frac{\log p}{T})^{-q}$, this bound is clearly sub-optimal when $T>\log p$ since the upper bound provided in our result scales as $R_q(\frac{\log p}{T})^{1-\frac{q}{2}}$.

\section{Simulation results}
\label{sc:simulation}
In order to assess our lower and upper bound results further, especially the dependence on $T$ and $n$, we conduct some simple numerical experiments.

The true signal $\theta^*$ is constructed as follows: (i) we set the first coefficient of $\theta^*$ to be $\theta^*_1=\frac{1}{\sqrt{p}}$ as we discussed earlier; (ii) since we assume that $D\theta^*\in\mathbb{R}^{p}_{+}$, then by following steps for constructing the packing set of $\mathcal{F}_{p,q,D}$ (described in Section \ref{ProofLower})   the remaining $p-1$ coefficients for $\theta^*$ are randomly generated from $\text{Unif}(0,\frac{1}{p\lambda_p})$, where $\lambda_p$ is the $p$-sparse localization quantity of $D$; (iii) the last step is to normalize $\bar{\theta}^*=[\theta_2^*,...,\theta_p^*]$ such that $\|\bar{\theta}^*\|_q^q\leq R_q.$ By following these steps we can generate signal $\theta^*$ such that $D\theta^*\in\mathcal{F}_{p,q,D}$. We choose two different basis (DCT and DWT) for $D$ in numerical experiments. Sensing matrix $A$ is generated using (\ref{eq:construction}) with re-scaled Bernoulli ensemble matrix $\tilde{A}$.

We compare the simple standard Lasso estimator analyzed in this paper (which we call \textbf{Lasso}) to two other comparison methods. Firstly we compare to the weighted Lasso estimator using data-dependent weights, that means we modify the estimator in (\ref{equation:l1problem}) by changing $\|\theta\|_1$ to $\sum_{i=1}^pd_k|\theta_k|$, the weights $d_k~(1\leq k\leq p)$ are chosen by using data-dependent technique stated in (\cite{Jiang2015data}). We call this estimator \textbf{WLasso}. Finally we consider the penalized Poisson likelihood estimator with $\ell_1$ regularization using SPIRAL algorithm developed in (\cite{harmany2012spiral}) which we call \textbf{PoissonLike}.

\begin{figure}
	\centering     
	\subfloat[DCT ($p = 1024,~R_q = 7,~q = 0.5,~T = 10^8$)]{\label{fig:an}\includegraphics[width=70mm]{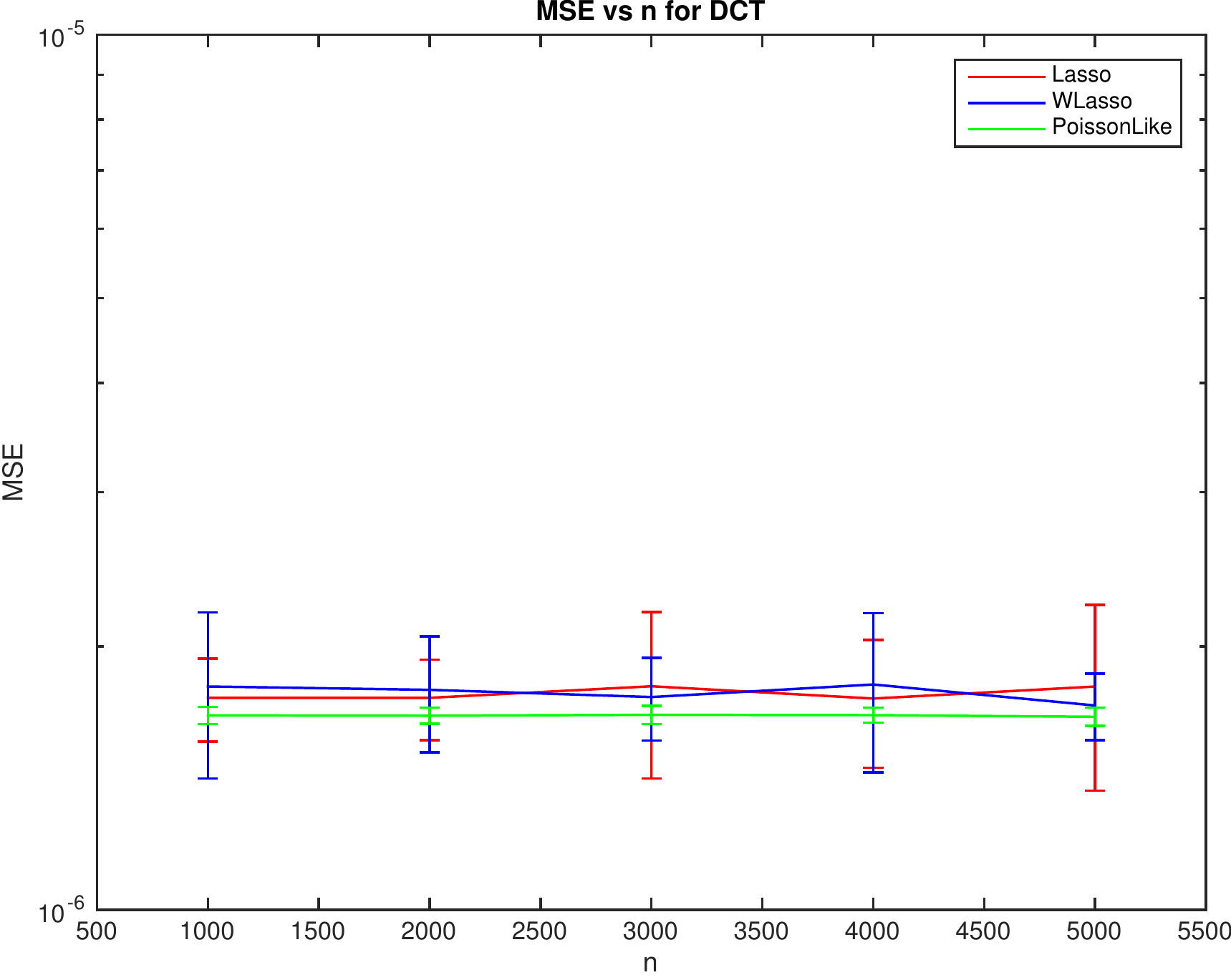}}
	\subfloat[DWT ($p = 1024,~R_q = 7,~q = 0.5,~T = 10^8$)]{\label{fig:bn}\includegraphics[width=70mm]{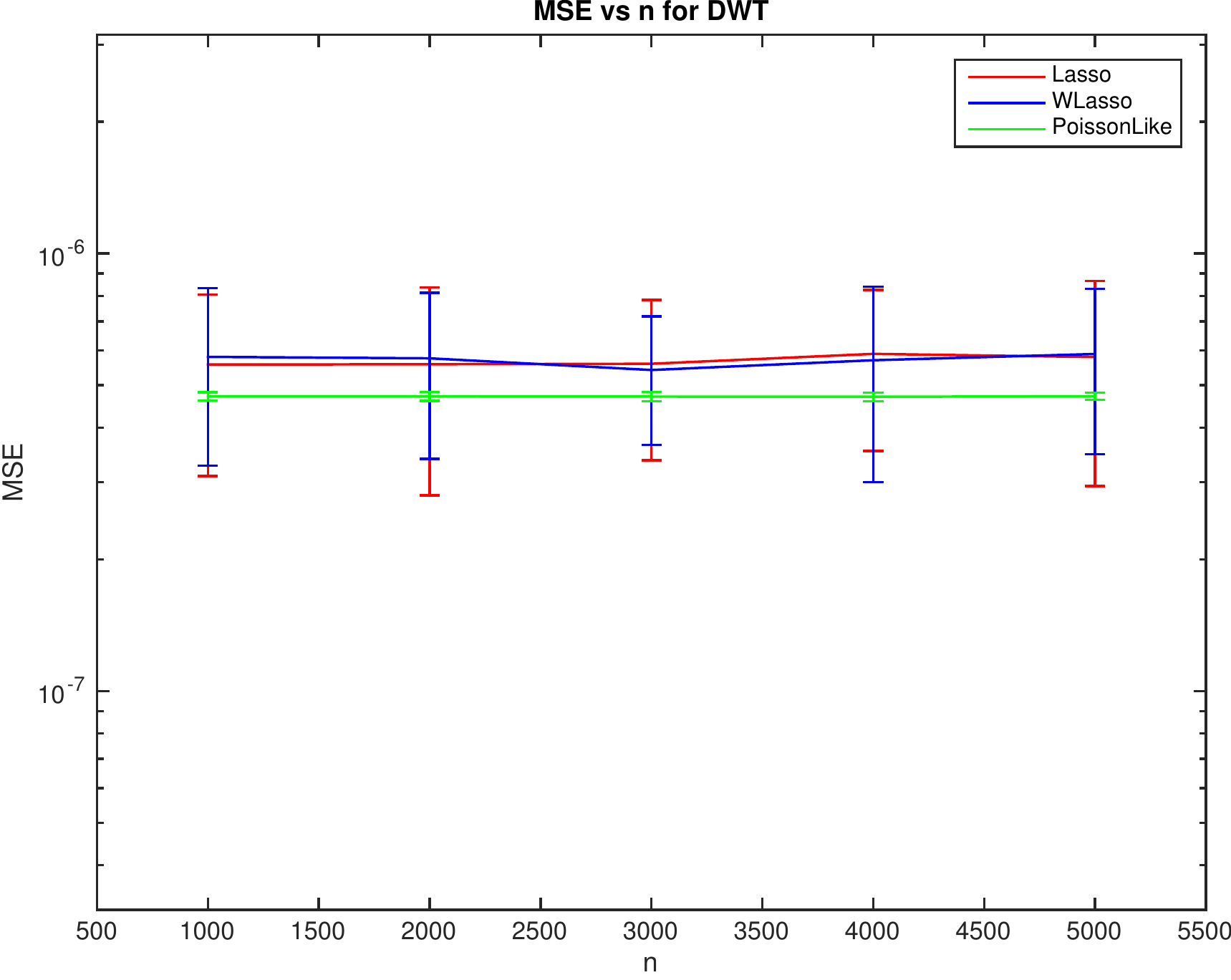}}
	\caption{MSE {\em vs.} $n$ for DCT and DWT}
	\label{figure:n}
\end{figure}

\begin{figure}
	\centering     
	\subfloat[DCT ($p = 1024,~R_q = 7,~q = 0.5,~n = 1000$)]{\label{fig:aT}\includegraphics[width=70mm]{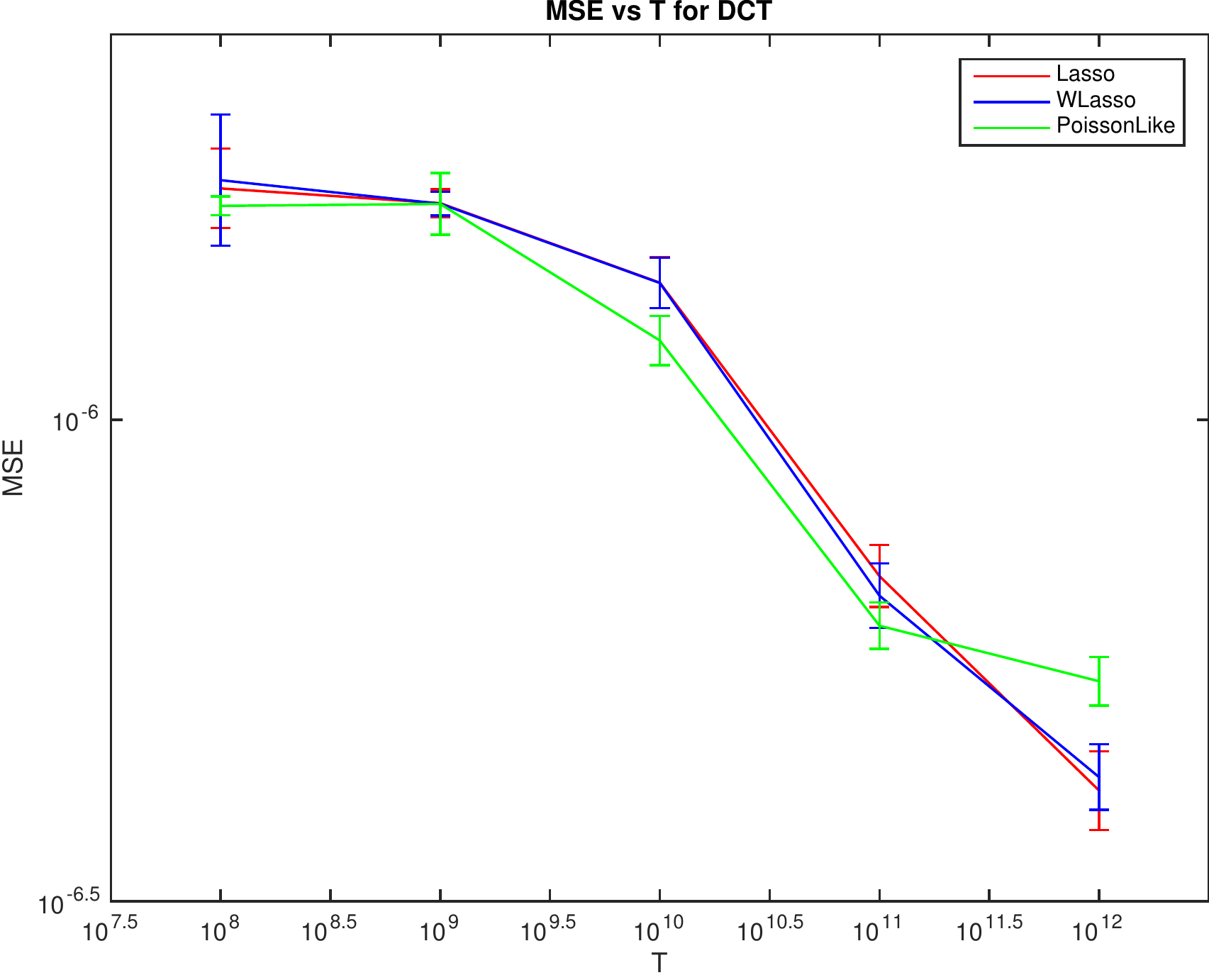}}
	\subfloat[DWT ($p = 1024,~R_q = 7,~q = 0.5,~n = 1000$)]{\label{fig:bT}\includegraphics[width=70mm]{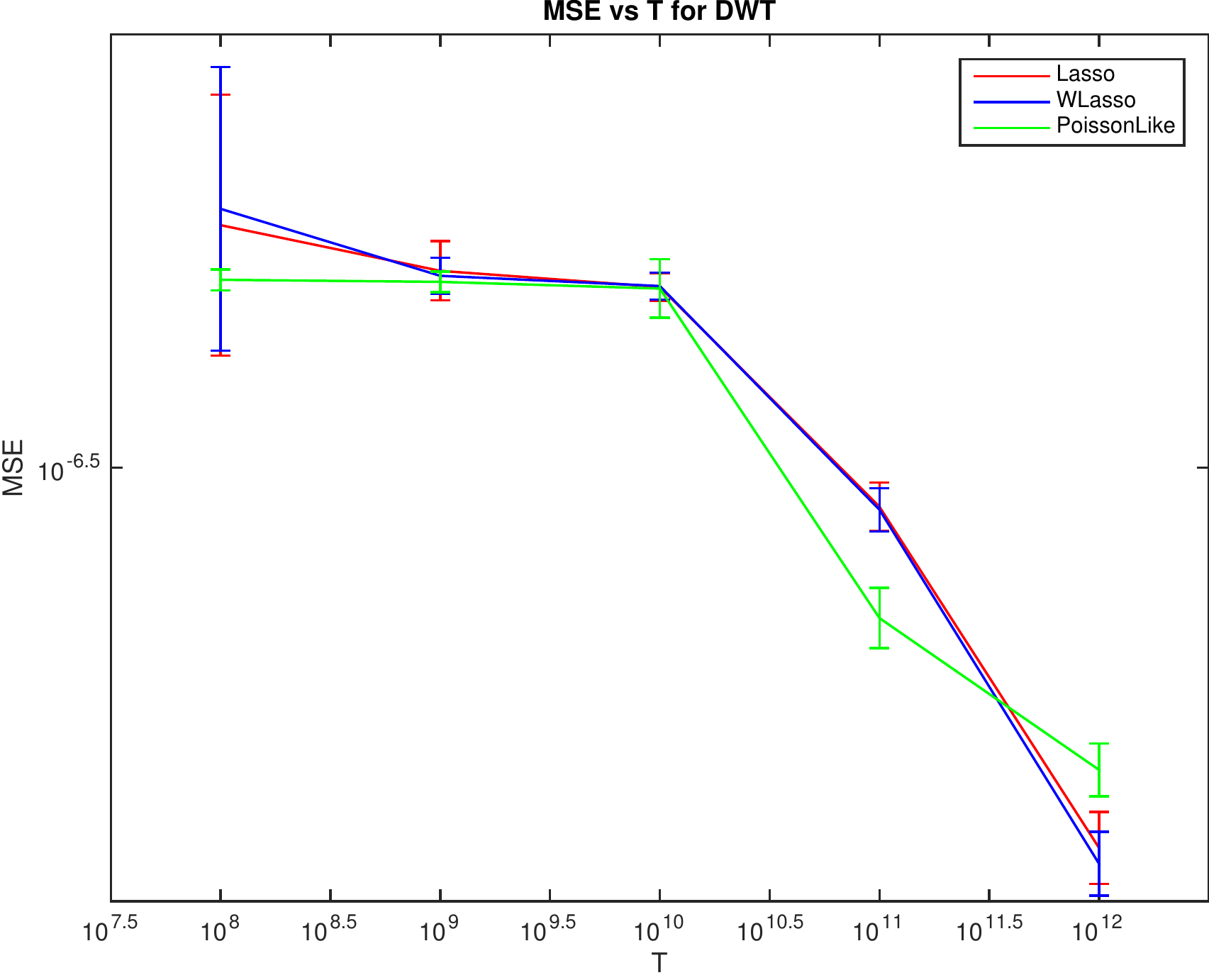}}
	\caption{MSE {\em vs.} $T$ for DCT and DWT}
	\label{figure:T}
\end{figure}

\begin{figure}
	\centering     
	\subfloat[DCT ($p = 1024,~R_q = 7,~n = 1000,~T = 10^8$)]{\label{fig:aq}\includegraphics[width=70mm]{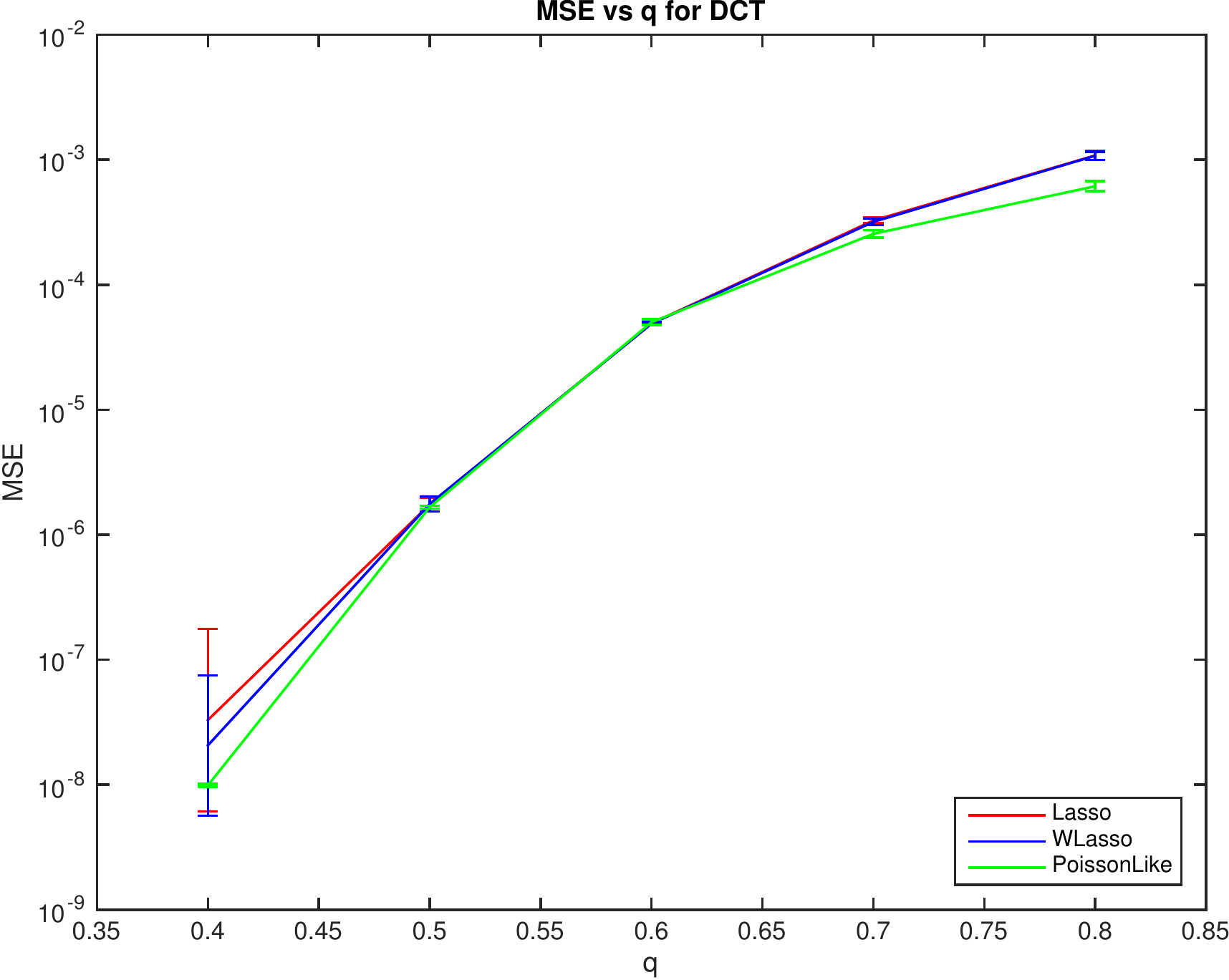}}
	\subfloat[DWT ($p = 1024,~R_q = 7,~n = 1000,~T = 10^8$)]{\label{fig:bq}\includegraphics[width=70mm]{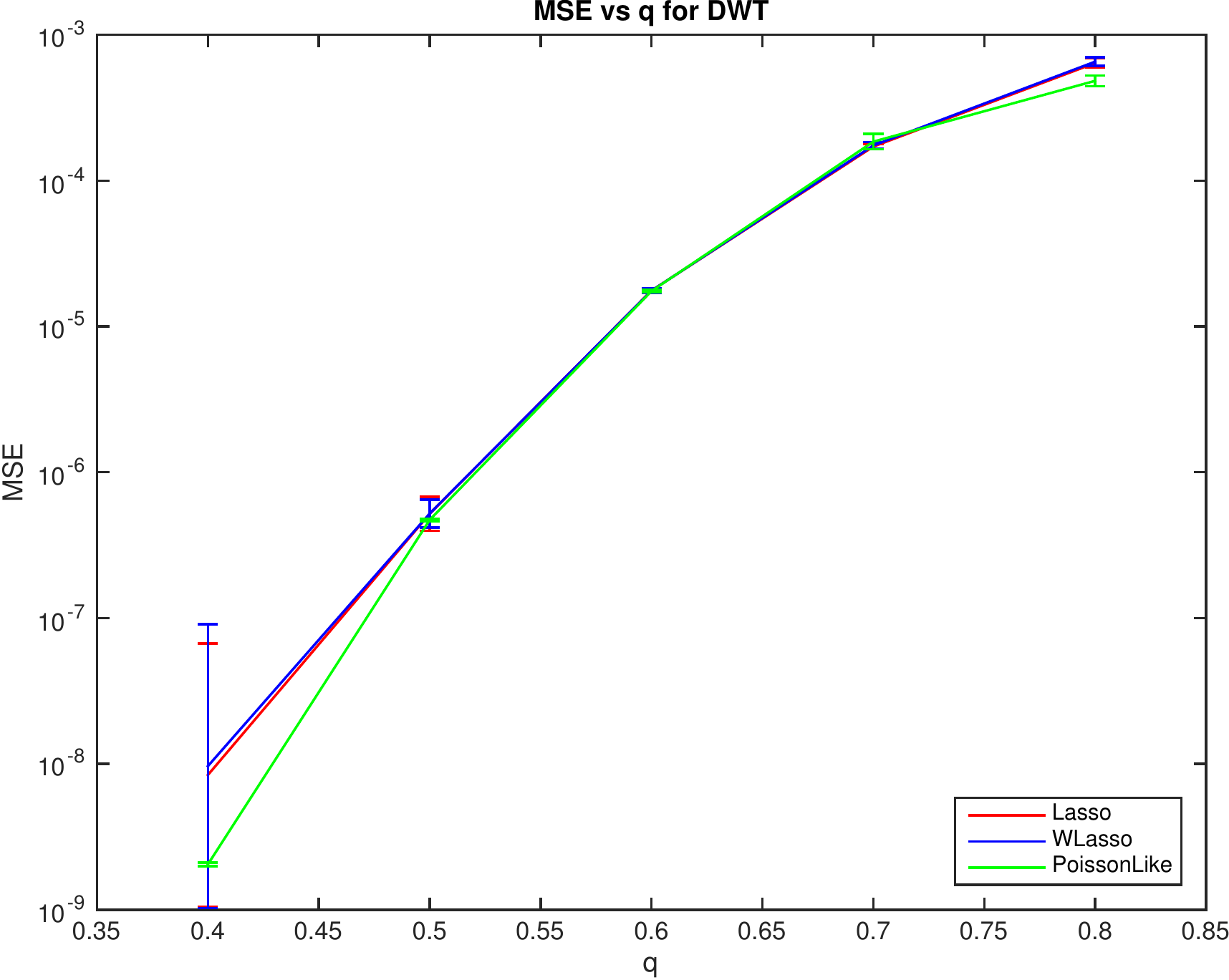}}
	\caption{MSE {\em vs.} $q$ for DCT and DWT}
	\label{figure:q}
\end{figure}

The results are shown from Figure \ref{figure:n} to Figure \ref{figure:q}. All tuning parameters are chosen by 5-fold cross validation. In these figures we show the average mean-squared error (MSE) results over 100 experiments for each data point. We can see all three methods have similar performance in terms of MSE. From Figure \ref{figure:n} we can see that the MSE results are independent of $n$ as long as $n$ is large enough to ensure that the sensing matrix $A$ satisfies the imposed assumptions. Figure \ref{figure:T} shows that MSE decreases with $T$ for $T$ large enough, this is consistent with our conclusion that the error bound results depend explicitly on the the intensity $T$ and not on the sample size $n$. Similar numerical results about the relationship of MSE with $T$ and $n$ can be seen in \cite{jiang2014minimax}. Figure \ref{figure:q} shows that MSE increases when $q$ increases, this is also reflected in our error bound results.

Note that the goal of these simulation results is to verify the theoretical results. Clearly all three methods perform well but to the best of our knowledge, the results in this paper are the first to provide theoretical guarantees for any estimator. It remains an open question to provide sharp theoretical guarantees for \textbf{WLasso} and \textbf{PoissonLike}. It also must be pointed out that we consider $\theta^*$ and $A$ that satisfy the assumptions specified. If these assumptions are violated, none of the three methods are guaranteed to perform well.
 
\section{Proofs}
\label{SecProofs}
In this section we provide the proofs for our three main results. We defer the more technical steps to the appendix.

\subsection{Proof of Theorem \ref{theorem:lower}}
\label{ProofLower}
The proof for the lower bound uses a combination of standard information-theoretic techniques involving Fano's inequality, and the explicit construction of a packing set that satisfies the $\ell_q$-ball constraint and our other physical constraints. In particular, the proof involves constructing a packing set for $\mathcal{F}_{p,q,D}$ and then applying the generalized Fano method to the packing set (see \cite{han1994generalizing}, \cite{IbrHas81} and \cite{YanBar99} for details). Constructing the packing set is the main challenge and novelty in the proof. Our packing set is based on a constrained hypercube construction in \cite{jiang2014minimax} along with the hyper-cube construction for $\ell_q$-balls in \cite{Kuh01}.

\begin{proof}
We begin our proof by constructing a packing set for $\mathcal{F}_{p,q,D}.$ For $1\leq k\leq \tilde{K}$ let $\mathcal{H}_k=\{\beta\in\{-1,0,+1\}^{p-1}:\|\beta\|_0=k\}.$ From [\cite{RasWaiYu11}, Lemma 4] we can find a subset $\tilde{\mathcal{H}}_k\subseteq\mathcal{H}_k$ with cardinality $|\tilde{\mathcal{H}}_k|\geq\text{exp}(\frac{k}{2}\text{log}\frac{p-\frac{k}{2}-1}{k})$ such that the Hamming distance $\rho_{H}(\beta,\beta')\geq\frac{k}{2}$ for all $\beta,\beta'\in\tilde{\mathcal{H}}_k$. Then note that $(\frac{R_q}{k})^{\frac{1}{q}}\mathcal{H}_k=\{\beta\in\{-(\frac{R_q}{k})^{\frac{1}{q}},0,+(\frac{R_q}{k})^{\frac{1}{q}}\}^{p-1}:\|\beta\|_q^q=R_q\}$. Now consider the re-scaled hypercube $\tilde{\mathcal{H}}_k$ by $\alpha_k$ with $0<\alpha_k\leq(\frac{R_q}{k})^{\frac{1}{q}}$, we define: $\mathcal{H}_{k,\alpha_k}=\{\theta\in\mathbb{R}^p:\theta=[1/\sqrt{p},\alpha_k\beta^{\top}]^{\top},\beta\in\tilde{\mathcal{H}}_k\}$, let $\eta^2_{\alpha_k}=\frac{k}{2}\alpha_k^2,$ then $\mathcal{H}_{k,\alpha_k}$ is a $\eta_{\alpha_k}$-packing set for $\mathcal{F}_{p,q,D}$ in the $\ell_2$ norm. To contrast with the packing set used in \cite{jiang2014minimax}, we require the extra constraint that each vertex has value at most $(\frac{R_q}{k})^{\frac{1}{q}}$ to ensure $\|\beta\|_q^q \leq R_q$.

The following lemma shows useful properties of this packing set.
\blems[\cite{jiang2014minimax}, Lemma 4.1]
\label{lemma:3prop}
For $1\leq k\leq \tilde{K}$, let $\lambda_k=\lambda_k(\bar{D}).$ Then the packing sets $\mathcal{H}_{k,\alpha_k}$ with $0<\alpha_k\leq\frac{1}{p\lambda_k}$ have the following properties:\\ \\
1. The $\ell_2$ distance between any two points $\theta$ and $\theta'$ in $\mathcal{H}_{k,\alpha_k}$ is bounded:
$$\eta_{\alpha_k}^2\leq\|\theta-\theta'\|^2_2\leq8\eta_{\alpha_k}^2.$$
2. For any $\theta\in\mathcal{H}_{k,\alpha_k}$, the corresponding $f=D\theta$~satisfies:
$$f_i\geq0,~\forall i\in\{1,2,...,p\},~~\text{and}~~~\|f\|_1=1.$$
3. The size of the packing set
$$|\mathcal{H}_{k,\alpha_k}|\geq\exp(\frac{k}{2}\log\frac{p-\frac{k}{2}-1}{k}).$$
\elems
The proof for this lemma can be found in \cite{jiang2014minimax}. 

For convenience we define the matrix $\Phi\triangleq AD$ and then $\Phi\theta=Af~\text{since}~f=D\theta$. Next we will apply the generalized Fano's method to the packing set, these techniques are developed in \cite{han1994generalizing}, \cite{IbrHas81} and \cite{YanBar99}. Define $M_k$ to be the cardinality of the set $\mathcal{H}_{k,\alpha_k}$, and the elements in $\mathcal{H}_{k,\alpha_k}$ can be denoted as $\{\theta^1,...,\theta^{M_k}\}$. Let $\tilde{\Theta}\in\mathbb{R}^p$ be a random vector drawn from a uniform distribution over the packing set $\{\theta^1,...,\theta^{M_k}\}$. Further let $\tilde{\theta} =\arg\min_{\theta \in \mathcal{H}_{k,\alpha_k}} \|\theta - D^{\top} f\|_2$. Then since $D$ is an orthonormal basis we can bound the minimax estimation error according to \cite{YanBar99}:
\begin{equation}
\label{eq:fano1}
P(\min_{f}\max_{f^{*}\in\mathcal{F}_{p,q,D}}\|f-f^{*}\|_2^2\geq\frac{\eta_{\alpha_k}^2}{4})\geq\min_{\tilde{\theta}}\mathbb{P}[\tilde{\theta}\neq\tilde{\Theta}].
\end{equation}
Applying Fano's inequality yields the following lower bound:
\begin{equation}
\label{eq:fano2}
\mathbb{P}[\tilde{\theta}\neq\tilde{\Theta}]\geq 1-\frac{I(y;\tilde{\Theta})+\log2}{\log M_k},
\end{equation}
where $y|\tilde{\Theta} \sim\text{Poisson}(T\Phi\tilde{\Theta})$ and $I(y;\tilde{\Theta})$ is the mutual information between random variable $y$ and $\tilde{\Theta}.$ Then from \cite{han1994generalizing} we have
\begin{equation}
\label{eq:fano3}
I(y;\tilde{\Theta})\leq\frac{1}{\binom{M_k}{2}}\sum_{i,j=1,..,M_k,i\neq j}\text{KL}(p(y|T\Phi\theta^{i})\|p(y|T\Phi\theta^{j})),
\end{equation}
where $\text{KL}(p_1\|p_2)$ is the Kullback-Leibler (KL) divergence between distributions $p_1$ and $p_2$. We will use the following lemma to bound KL divergence of Poisson distributions in terms of the squared $\ell_2$-distance.
\blems[\cite{jiang2014minimax}, Lemma 4.2]
\label{Lemma:KL}
Let $p(y|\mu)$ denote the vector Poisson distribution with mean parameter $\mu\in\mathbb{R}^{n}_{+}.$~For~$\mu_1,\mu_2\in\mathbb{R}^{n}_{+}$, if there exists some value $c>0$ such that $\mu_2\succeq c{\mathbbm{1}}_{n\times1}$, then the following holds:
$$\text{KL}(p(y|\mu_1)\|p(y|\mu_2))\leq\frac{1}{c}\|\mu_1-\mu_2\|^2_2.$$
\elems

The following lemma shows that entries in $Af^*$ are bounded between $\frac{1}{2n}$ and $\frac{1}{n}$ under Assumption \ref{as:construction}:
\blems[\cite{jiang2014minimax}, Lemma 4.3]
\label{lemma:Afbound}
If the sensing matrix A satisfies Assumption \ref{as:construction}, then for all non-negative $f$ with $\|f\|_1=1$, we have:
$$\frac{1}{2n}\mathbbm{1}_{n\times1}\preceq Af \preceq\frac{1}{n}\mathbbm{1}_{n\times1}.$$
\elems

The proofs for Lemmas~\ref{Lemma:KL} and ~\ref{lemma:Afbound} can be found in \cite{jiang2014minimax}. By Lemma \ref{lemma:Afbound} we have $\Phi\theta^{j}=AD\theta^j\succeq\frac{1}{2n}\mathbbm{1}_{n\times1}$ and then it follows that $T\Phi\theta^{j}\succeq\frac{T}{2n}\mathbbm{1}_{n\times1}$. Then from Lemma \ref{Lemma:KL} we can bound the KL divergence between $p(y|T\Phi\theta^{i})$ and $p(y|T\Phi\theta^{j})$ as follows:
\begin{equation}
\label{eq:fano4}
\text{KL}(p(y|T\Phi\theta^{i})\|p(y|T\Phi\theta^{j}))\leq\frac{2n}{T}\|T\Phi(\theta^{i}-\theta^{j})\|^2_2=2nT\|\Phi(\theta^{i}-\theta^{j})\|^2_2.
\end{equation}

By Assumption \ref{as:construction} and \ref{as:rip} if we denote $f^{i}=D\theta^{i},~f^{j}=D\theta^{j}$, then 
\begin{eqnarray}
\|\Phi(\theta^{i}-\theta^{j})\|^2_2 & = & \|A(f^{i}-f^{j})\|^2_2 \nonumber \\ 
& = &  \|\frac{1}{2(a_u-a_\ell)\sqrt{n}}\tilde{A}D(\theta^{i}-\theta^{j})\|^2_2 \nonumber \\
& \leq & \frac{1+\delta_{\tilde{K}}}{4(a_u-a_\ell)^2n}\|\theta^{i}-\theta^{j}\|_2^2.
\end{eqnarray}
Since $\|\theta^i-\theta^j\|^2_2\leq 8\eta_{\alpha_k}^2$ by Lemma \ref{lemma:3prop}, we further have
\begin{equation}
\label{eq:fano5}
\|\Phi(\theta^{i}-\theta^{j})\|^2_2\leq\frac{2(1+\delta_{\tilde{K}})}{(a_u-a_\ell)^2n}\eta_{\alpha_k}^2.
\end{equation}

Then by combining (\ref{eq:fano4}) and (\ref{eq:fano5}) we have:
\begin{equation}
\label{eq:fano6}
\text{KL}(p(y|T\Phi\theta^{i})\|p(y|T\Phi\theta^{j}))\leq 2nT\|\Phi(\theta^{i}-\theta^{j})\|^2_2\leq\frac{4(1+\delta_{\tilde{K}})T}{(a_u-a_\ell)^2}\eta_{\alpha_k}^2.
\end{equation}
Then the mutual information can be bounded by using (\ref{eq:fano3}) and (\ref{eq:fano6})
\begin{eqnarray}
I(y;\tilde{\Theta}) & \leq & \frac{1}{\binom{M_k}{2}}\sum_{i\neq j}\text{KL}(p(y|T\Phi\theta^{i})\|p(y|T\Phi\theta^{j}))\nonumber \\ & \leq & \max_{i\neq j}\text{KL}(p(y|T\Phi\theta^{i})\|p(y|T\Phi\theta^{j}))\nonumber \\
& \leq & \frac{4(1+\delta_{\tilde{K}})T}{(a_u-a_\ell)^2}\eta_{\alpha_k}^2.
\end{eqnarray}
Using (\ref{eq:fano2}) and the lower bound for $M_k$ we have
\begin{equation}
\label{eq:fano7}
\mathbb{P}[\tilde{\theta}\neq\tilde{\Theta}]\geq1-\frac{I(y;\tilde{\Theta})+\log2}{\log M_k}\geq1-\frac{\frac{4(1+\delta_{\tilde{K}})T}{(a_u-a_\ell)^2}\eta_{\alpha_k}^2+\log2}{\frac{k}{2}\log\frac{p-\frac{k}{2}-1}{k}}.
\end{equation}
Next we will show that the probability in (\ref{eq:fano7}) is bounded by the constant $1/2$. This constant is guaranteed if the following two inequalities are true:

\begin{eqnarray}
\label{equation:kineq1}
\frac{k}{2}\log\frac{p-\frac{k}{2}-1}{k} & \geq & 4\log2,\\
\label{equation:kineq2}
\frac{k}{2}\log\frac{p-\frac{k}{2}-1}{k} & \geq & \frac{16(1+\delta_{\tilde{K}})T}{(a_u-a_\ell)^2}\eta_{\alpha_k}^2.
\end{eqnarray}
For the first inequality (\ref{equation:kineq1}) if $k=1$,
$$\frac{k}{2}\log\frac{p-\frac{k}{2}-1}{k}=\frac{1}{2}\log{(p-\frac{3}{2})}\geq\frac{1}{2}\log{(\frac{515}{2}-\frac{3}{2})}=4\log2,$$
where the inequality is a result of $p\geq260$. And if $k\geq2$,
\begin{eqnarray*}
\frac{k}{2}\log\frac{p-\frac{k}{2}-1}{k} & \geq & \log\frac{p-\frac{k}{2}-1}{k}\\
& \geq & \log\frac{p-\frac{\tilde{K}}{2}-1}{\tilde{K}}\\
& \geq & \log\frac{(\frac{33\tilde{K}}{2}+1)-(\frac{\tilde{K}}{2}+1)}{\tilde{K}}\\
& = & 4\log2,
\end{eqnarray*}
where the inequality is the result of $p\geq\frac{33\tilde{K}}{2}+1$. For the second inequality (\ref{equation:kineq2}) we need:
$$\frac{k}{32}\log\frac{p-\frac{k}{2}-1}{k}\geq\frac{(1+\delta_{\tilde{K}})T}{(a_u-a_\ell)^2}\eta_{\alpha_k}^2,$$
which leads to
$$\eta_{\alpha_k}^2\leq\frac{(a_u-a_\ell)^2k}{32(1+\delta_{\tilde{K}})T}\log\frac{p-\frac{k}{2}-1}{k}.$$
Since for Lemma \ref{lemma:3prop} we require that $0<\alpha_k\leq\frac{1}{p\lambda_k}$, also recall the extra constraint that $\alpha_k\leq(\frac{R_q}{k})^{\frac{1}{q}}$, thus we have:
$$\eta_{\alpha_k}^2=\min\left(\frac{k}{2}(\frac{1}{p\lambda_k})^2,\frac{(a_u-a_\ell)^2k}{32(1+\delta_{\tilde{K}})T}\log\frac{p-\frac{k}{2}-1}{k},\frac{k}{2}(\frac{R_q}{k})^{\frac{2}{q}}\right).$$

Then with probability greater than $\frac{1}{2}$ we have
\begin{eqnarray}
\label{equation:lowerbound}
\min_{f}\max_{f^{*}\in\mathcal{F}_{p,q,D}}\|f-f^{*}\|_2^2 & \geq & \frac{\eta_{\alpha_k}^2}{4} \nonumber \\
& = & \min\left(\frac{k}{8}(\frac{1}{p\lambda_k})^2,\frac{(a_u-a_\ell)^2k}{128(1+\delta_{\tilde{K}})T}\log\frac{p-\frac{k}{2}-1}{k},\frac{k}{8}(\frac{R_q}{k})^{\frac{2}{q}}\right).
\end{eqnarray}

In order to further simplify (\ref{equation:lowerbound}), note that $k\leq\tilde{K}=O(R_q(\frac{\log p}{T})^{-\frac{q}{2}})$, we then have 
\begin{equation}
\label{equation:tildeK}
k\leq O(R_q(\frac{\log\frac{p}{k}}{T})^{-\frac{q}{2}}),
\end{equation}
from (\ref{equation:tildeK}) there exists some absolute constant $C_1>0$ such that
\begin{equation}
\label{equation:lasttwo}
\frac{(a_u-a_\ell)^2k}{128(1+\delta_{\tilde{K}})T}\log\frac{p-\frac{k}{2}-1}{k}\leq C_1\frac{k}{8}(\frac{R_q}{k})^{\frac{2}{q}}.
\end{equation}
Thus by (\ref{equation:lowerbound}) and (\ref{equation:lasttwo})
\begin{equation}
\label{equation:lowerbound2}
\min_{f}\max_{f^{*}\in\mathcal{F}_{p,q,D}}\|f-f^{*}\|_2^2\geq\min\left(\frac{k}{8}(\frac{1}{p\lambda_k})^2,\frac{(a_u-a_\ell)^2k}{128\max\{1,C_1\}(1+\delta_{\tilde{K}})T}\log\frac{p-\frac{k}{2}-1}{k}\right).
\end{equation}
Since (\ref{equation:lowerbound2}) holds for all $1\leq k\leq\tilde{K}$, there exists a constant $C_L>0$ such that
$$\min_{f}\max_{f^{*}\in\mathcal{F}_{p,q,D}}\|f-f^{*}\|_2^2\geq C_L \max_{1 \leq k \leq \tilde{K}}\left\{\min\left(\frac{k}{p^2\lambda_{k}^2},\frac{k\log\frac{p}{k}}{T}\right)\right \}$$
with probability greater than $\frac{1}{2}.$
\end{proof}

\subsection{Proof for Theorem \ref{ThmL1}}
\label{SecProofL1}
The proof for the upper bound involves direct analysis of the lasso estimator defined in (\ref{equation:l1problem}). Our analysis follows standard steps for analysis of regularized M-estimators (see \cite{BiRiTsy08,Neg10,vandeGeer}) along with addressing two challenges specific to this setting:
(1) we use concentration bounds for linear combination of Poisson random variables and how they are used to determine a $\lambda_n$; (2) use Assumption \ref{Assumption:ntildek} to show that matrix $A\bar{D}$ satisfies the restricted eigenvalue condition and satisfies the physical constraints. 

\spro
From (\ref{equation:l1problem}) in Section \ref{Secl1method} we know $\hat{\theta}_{\lambda_n}$ is a solution to the following problem:
\begin{equation}
\label{eq:l1problem}
\hat{\theta}_{\lambda_n} \in \text{arg}\min_{\theta\in\mathbb{R}^p,\theta_1=\frac{1}{\sqrt{p}}}\frac{n}{T^2}\|y-TAD\theta\|_2^2+\lambda_n\|\theta\|_1,
\end{equation}
 with $(\hat{\theta}_{\lambda_n})_1=\frac{1}{\sqrt{p}}.$ Since $\theta^*$ satisfies the constraint that $\theta^*\in\mathbb{R}^p$ and $\theta^*_1=\frac{1}{\sqrt{p}},$ we have the following basic inequality
$$\frac{n}{T^2}\|y-TAD\hat{\theta}_{\lambda_n}\|_2^2+\lambda_n\|\hat{\theta}_{\lambda_n}\|_1\leq\frac{n}{T^2}\|y-TAD\theta^*\|_2^2+\lambda_n\|\theta^*\|_1.$$
Hence
$$\frac{n}{T^2}\|TAD(\theta^*-\hat{\theta}_{\lambda_n})\|_2^2\leq\frac{2n}{T^2}|(y-TAD\theta^*)^{\top}TAD(\theta^*-\hat{\theta}_{\lambda_n})|+\lambda_n(\|\theta^*\|_1-\|\hat{\theta}_{\lambda_n}\|_1),$$
and then
\begin{equation}
\label{equation:lassodiff1}
n\|AD(\theta^*-\hat{\theta}_{\lambda_n})\|_2^2\leq\frac{2n}{T}|(y-TAD\theta^*)^{\top}AD(\theta^*-\hat{\theta}_{\lambda_n})|+\lambda_n(\|\theta^*\|_1-\|\hat{\theta}_{\lambda_n}\|_1).
\end{equation}
Note that $(\theta^*)_1=(\hat{\theta}_{\lambda_n})_1=\frac{1}{\sqrt{p}}$, thus it is reasonable to define the error vector $\hat{\Delta}=\bar{\theta}^*-\bar{\hat{\theta}}_{\lambda_n}\in\mathbb{R}^{p-1}$, where $\bar{\theta}^*=[\theta^*_2,...,\theta^*_p]^{\top}\in\mathbb{R}^{p-1}$ and $\bar{\hat{\theta}}_{\lambda_n}=[(\hat{\theta}_{\lambda_n})_2,...,(\hat{\theta}_{\lambda_n})_p]^{\top}\in\mathbb{R}^{p-1}$. Then (\ref{equation:lassodiff1}) can be reduced to:
\begin{eqnarray}
n\|A\bar{D}\hat{\Delta}\|_2^2 & \leq & \frac{2n}{T}|(y-TAD\theta^*)^{\top}A\bar{D}\hat{\Delta}|+\lambda_n(\|\theta^*\|_1-\|\hat{\theta}_{\lambda_n}\|_1) \nonumber\\
\label{equation:lassodiff2}
& \leq & \|\frac{2n}{T}(y-TAD\theta^*)^{\top}A\bar{D}\|_{\infty}\|\hat{\Delta}\|_1+\lambda_n(\|\theta^*\|_1-\|\hat{\theta}_{\lambda_n}\|_1),
\end{eqnarray}
where $\bar{D}=[d_2,...,d_p]\in\mathbb{R}^{p\times(p-1)}$.

In order to associate the term $\|\theta^*\|_1-\|\hat{\theta}_{\lambda_n}\|_1$ with $\hat{\Delta}$, we define a threshold parameter $\eta>0$ and the threshold subset as follows:
$$S_{\eta}:=\{j\in\{2,3,...,p\}~|~|\theta^*_j|>\eta\}$$
and its complement
$$S^c_{\eta}:=\{j\in\{2,3,...,p\}~|~|\theta^*_j|\leq\eta\}.$$
Suppose $u$ is a vector in $\mathbb{R}^{p-1}$, we will define $u_{S_\eta}\in\mathbb{R}^{p-1}$ as following:
$$(u_{S_\eta})_j=\begin{cases}
u_j~~\text{if}~j+1\in S_{\eta}\\
0~~\text{if}~j+1\not\in S_{\eta}
\end{cases}
~\text{for}~1\leq j\leq p-1,
$$
and $u_{S^c_{\eta}}$ is defined in a similar way.
Now we show how to connect $\|\theta^*\|_1-\|\hat{\theta}_{\lambda_n}\|_1$ with $\hat{\Delta}$. Note that 
\begin{equation}
\label{equation:proofdevia} 
\|\theta^*\|_1-\|\hat{\theta}_{\lambda_n}\|_1=\|\bar{\theta}^*\|_1-\|\bar{\hat{\theta}}_{\lambda_n}\|_1,
\end{equation}
since $\theta_1^* = (\hat{\theta}_{\lambda_n})_1 = \frac{1}{\sqrt{p}}$.
Then by using the triangle inequality we have
\begin{eqnarray}
\|\bar{\hat{\theta}}_{\lambda_n}\|_1=\|\bar{\theta}^*-\hat{\Delta}\|_1 &=& \|\bar{\theta}^*_{S_{\eta}}+\bar{\theta}^*_{S^c_{\eta}}-\hat{\Delta}_{S_{\eta}}-\hat{\Delta}_{S^c_{\eta}}\|_1\nonumber \\
& \geq & \|\bar{\theta}^*_{S_{\eta}}-\hat{\Delta}_{S^c_{\eta}}\|_1-\|\bar{\theta}^*_{S^c_{\eta}}\|_1-\|\hat{\Delta}_{S_{\eta}}\|_1\nonumber\\
&=& \|\bar{\theta}^*_{S_{\eta}}\|_1+\|\hat{\Delta}_{S^c_{\eta}}\|_1-\|\bar{\theta}^*_{S^c_{\eta}}\|_1-\|\hat{\Delta}_{S_{\eta}}\|_1.
\end{eqnarray}
On the other hand we have $\|\bar{\theta}^*\|_1\leq\|\bar{\theta}^*_{S_{\eta}}\|_1+\|\bar{\theta}^*_{S^c_{\eta}}\|_1$.
Thus by combining these two inequalities we have
$$\|\bar{\theta}^*\|_1-\|\bar{\hat{\theta}}_{\lambda_n}\|_1\leq\|\hat{\Delta}_{S_{\eta}}\|_1-\|\hat{\Delta}_{S^c_{\eta}}\|_1+2\|\bar{\theta}^*_{S^c_{\eta}}\|_1.$$
Therefore by (\ref{equation:proofdevia}):
\begin{equation}
\label{equation:deviation}
\|\theta^*\|_1-\|\hat{\theta}_{\lambda_n}\|_1\leq\|\hat{\Delta}_{S_{\eta}}\|_1-\|\hat{\Delta}_{S^c_{\eta}}\|_1+2\|\bar{\theta}^*_{S^c_{\eta}}\|_1.
\end{equation}
By using (\ref{equation:deviation}) in (\ref{equation:lassodiff2}) we have
\begin{equation}
\label{equation:lassodiff3}
n\|A\bar{D}\hat{\Delta}\|_2^2\leq\|\frac{2n}{T}(y-TAD\theta^*)^{\top}A\bar{D}\|_{\infty}\|\hat{\Delta}\|_1+\lambda_n(\|\hat{\Delta}_{S_{\eta}}\|_1-\|\hat{\Delta}_{S^c_{\eta}}\|_1+2\|\bar{\theta}^*_{S^c_{\eta}}\|_1).
\end{equation}
Now we upper bound the $\|.\|_{\infty}$ norm through the following Lemma:
\blems
\label{lemma:lambdan}
Under Assumption \ref{as:construction} and \ref{as:rip}, with probability at least $1-\frac{2}{p-1}$ that 
$$\|\frac{2n}{T}(y-TAD\theta^*)^{\top}A\bar{D}\|_{\infty}\leq\sqrt{\frac{32M\log p}{T}},$$
where $M=\frac{1+\delta_{\tilde{K}}}{4(a_u-a_{\ell})^2}$.
\elems
The proof for Lemma \ref{lemma:lambdan} is deferred to the appendix.

Thus by setting $\lambda_n=2\sqrt{\frac{32M\log p}{T}}$ in (\ref{equation:lassodiff3}) we have
\begin{eqnarray}
n\|A\bar{D}\hat{\Delta}\|_2^2 &\leq& \frac{\lambda_n}{2}\|\hat{\Delta}\|_1+\lambda_n(\|\hat{\Delta}_{S_{\eta}}\|_1-\|\hat{\Delta}_{S^c_{\eta}}\|_1+2\|\bar{\theta}^*_{S^c_{\eta}}\|_1)\nonumber \\
&\leq& \frac{\lambda_n}{2}(\|\hat{\Delta}_{S_{\eta}}\|_1+\|\hat{\Delta}_{S^c_{\eta}}\|_1+2\|\hat{\Delta}_{S_{\eta}}\|_1-2\|\hat{\Delta}_{S^c_{\eta}}\|_1+4\|\bar{\theta}^*_{S^c_{\eta}}\|_1) \nonumber \\
\label{equation:lassodiff4}
&=& \frac{\lambda_n}{2}(3\|\hat{\Delta}_{S_{\eta}}\|_1-\|\hat{\Delta}_{S^c_{\eta}}\|_1+4\|\bar{\theta}^*_{S^c_{\eta}}\|_1),
\end{eqnarray}
where the second inequality follows from the triangle inequality. From (\ref{equation:lassodiff4}) we can see that $0\leq3\|\hat{\Delta}_{S_{\eta}}\|_1-\|\hat{\Delta}_{S^c_{\eta}}\|_1+4\|\bar{\theta}^*_{S^c_{\eta}}\|_1$, then the error vector $\hat{\Delta}$ should satisfy $\|\hat{\Delta}_{S^c_{\eta}}\|_1\leq3\|\hat{\Delta}_{S_{\eta}}\|_1+4\|\bar{\theta}^*_{S^c_{\eta}}\|_1$. The following lemma shows that $n\|A\bar{D}\hat{\Delta}\|_2^2$ is lower bounded for all $\hat{\Delta}\in\{\Delta\in\mathbb{R}^{p-1}~|~\|\Delta_{S^c_\eta}\|_1\leq3\|\Delta_{S_\eta}\|_1+4\|\bar{\theta}^*_{S^c_\eta}\|_1\}$:
\blems
\label{lemma:rsc}
Suppose Assumption \ref{Assumption:ntildek} and \ref{Assumption:L1RE} hold, for all $\hat{\Delta}\in\{\Delta\in\mathbb{R}^{p-1}~|~\|\Delta_{S^c_\eta}\|_1\leq 3\|\Delta_{S_\eta}\|_1+4\|\bar{\theta}^*_{S^c_\eta}\|_1\}$ we have
$$n\|A\bar{D}\hat{\Delta}\|^2_2\geq c_1k_1^2\|\hat{\Delta}\|^2_2-c_2k_2^2\frac{\log p}{n}\|\bar{\theta}^*_{S^c_\eta}\|_1^2,$$
where $c_1,c_2>0$ are some constants. 
\elems
Once again the proof of Lemma \ref{lemma:rsc} is deferred to the appendix.

By using Lemma \ref{lemma:rsc} in (\ref{equation:lassodiff4}) we have
\begin{eqnarray}
c_1k_1^2\|\hat{\Delta}\|^2_2-c_2k_2^2\frac{\log p}{n}\|\bar{\theta}^*_{S^c_\eta}\|_1^2 &\leq& \frac{\lambda_n}{2}(3\|\hat{\Delta}_{S_{\eta}}\|_1-\|\hat{\Delta}_{S^c_{\eta}}\|_1+4\|\bar{\theta}^*_{S^c_{\eta}}\|_1)\nonumber \\
\label{equation:lassodiff5}
&\leq& \frac{\lambda_n}{2}(3\|\hat{\Delta}_{S_{\eta}}\|_1+4\|\bar{\theta}^*_{S^c_{\eta}}\|_1).
\end{eqnarray}
Then note that $\|\hat{\Delta}_{S_{\eta}}\|_1\leq\sqrt{|S_\eta|}\|\hat{\Delta}_{S_{\eta}}\|_2\leq\sqrt{|S_\eta|}\|\hat{\Delta}\|_2$, where $|S_{\eta}|$ is the cardinality of set $S_{\eta},$ then from (\ref{equation:lassodiff5}) we have
$$c_1k_1^2\|\hat{\Delta}\|^2_2-c_2k_2^2\frac{\log p}{n}\|\bar{\theta}^*_{S^c_\eta}\|_1^2\leq\frac{\lambda_n}{2}(3\sqrt{|S_\eta|}\|\hat{\Delta}\|_2+4\|\bar{\theta}^*_{S^c_{\eta}}\|_1),$$
which implies that
\begin{equation}
\label{equation:lassodiff6}
c_1k_1^2\|\hat{\Delta}\|^2_2-\frac{3\lambda_n\sqrt{|S_{\eta}|}}{2}\|\hat{\Delta}\|_2-c_2k_2^2\frac{\log p}{n}\|\bar{\theta}^*_{S^c_\eta}\|_1^2-2\lambda_n\|\bar{\theta}^*_{S^c_{\eta}}\|_1\leq0.
\end{equation}
Note that the left hand side of (\ref{equation:lassodiff6}) can be seen as a quadratic form of $\|\hat{\Delta}\|_2$. Thus by solving this quadratic inequality for $\|\hat{\Delta}\|_2$ we have
$$\|\hat{\Delta}\|_2^2\leq 9\frac{\lambda^2_n}{c_1^2k_1^4}|S_{\eta}|+\frac{1}{c_1k^2_1}(2c_2k_2^2\frac{\log p}{n}\|\bar{\theta}^*_{S^c_{\eta}}\|_1^2+4\lambda_n\|\bar{\theta}^*_{S^c_{\eta}}\|_1).$$
Hence
\begin{eqnarray}
\|\hat{f}_{\lambda_n}-f^*\|_2^2 &=& \|D(\hat{\theta}_{\lambda_n}-\theta^*)\|_2^2=\|\hat{\theta}_{\lambda_n}-\theta^*\|^2_2\nonumber\\
&=& \|\bar{\hat{\theta}}_{\lambda_n}-\bar{\theta}^*\|_2^2=\|\hat{\Delta}\|_2^2\nonumber \\
\label{equation:fbound}
&\leq& 9\frac{\lambda_n^2}{c_1^2k_1^4}|S_{\eta}|+\frac{1}{c_1k_1^2}(2c_2k_2^2\frac{\log p}{n}\|\bar{\theta}^*_{S^c_{\eta}}\|_1^2+4\lambda_n\|\bar{\theta}^*_{S^c_{\eta}}\|_1).
\end{eqnarray}
Since
\begin{equation}
\label{equation:Seta}
R_q\geq\sum_{j\in S^c_{\eta}}|{\theta}^*_j|^q\geq\eta^q|S_\eta|,
\end{equation}
we have~$|S_\eta|\leq\eta^{-q}R_q$. On the other hand
\begin{equation}
\label{equation:theta}
\|\bar{\theta}^*_{S^c_{\eta}}\|_1=\sum_{j\in S^c_{\eta}}|\bar{\theta}^*_j|=\sum_{j\in S^c_{\eta}}|\bar{\theta}^*_j|^q|\bar{\theta}^*_j|^{1-q}\leq R_q\eta^{1-q}.
\end{equation}
If we set $\eta=\frac{\lambda_n}{c_1k_1^2}$, by using (\ref{equation:Seta}) and (\ref{equation:theta}) in (\ref{equation:fbound}) we have
\begin{eqnarray}
\|\hat{f}_{\lambda_n}-f^*\|_2^2&\leq& 13R_q(\frac{\lambda_n}{c_1k_1^2})^{2-q}+\frac{2c_2k_2^2}{c_1k_1^2}\frac{\log p}{n}R_q^2(\frac{\lambda_n}{c_1k_1^2})^{2-2q}\\
\label{equation:fbound2}
&\leq&13R_q(\frac{\lambda_n}{c_1k_1^2})^{2-q}+\frac{2c_2k_2^2}{c_1k_1^2}R_q(\frac{\lambda_n}{c_1k_1^2})^{2-q}\left[\frac{R_q\log p}{n}\lambda_n^{-q}\right].
\end{eqnarray}

Since $\lambda_n=2\sqrt{\frac{32 M\log p}{T}}$, then under Assumption \ref{Assumption:ntildek} we can see that the term $\frac{R_q\log p}{n}\lambda_n^{-q}$ in (\ref{equation:fbound2}) scales as $O(1)$. Then with probability at least $1-\frac{2}{p-1}$ there exists a $C_U>0$ such that 
$$\|\hat{f}_{\lambda_n}-f^*\|_2^2\leq C_UR_q(\frac{\log p}{T})^{1-\frac{q}{2}}.$$
\fpro

\subsection{Proof for Theorem \ref{Theorem:RE}}
\label{SecProofRE}

The proof for Theorem \ref{Theorem:RE} uses techniques developed in \cite{RasWaiYu10b} adapted from Gaussian to  sub-Gaussian ensembles. The reason we adapt to sub-Gaussian ensembles is so that we construct a random ensemble that satisfies all the physical constraints. In the proof of \cite{RasWaiYu10b} the first step is to show the term $M(r,\Gamma):=\sup_{x\in V(r)}\{1-\frac{\|\Gamma x\|_2}{\sqrt{n}}\}$ is sharply concentrated around its expectation with high probability when $\Gamma$ is a matrix with Gaussian random variables, we will use [\cite{mendelson2007reconstruction}, Theorem 2.3] to show this is also true when $\Gamma$ is a matrix with subgaussian random variables. Finally we use peeling techniques to complete the proof, which are used in \cite{RasWaiYu10b}.

To begin we define the standard Gaussian width of a star-shaped set $T$. A set $T$ is \emph{start-shaped} if $cT \subset T$ for all $0 \leq c \leq 1$.
\bde[\cite{mendelson2007reconstruction}, Definition 2.1]
Let $T\subset\mathbb{R}^p$ and let $g_1,...,g_p$ be independent standard Gaussian random variables. Denote by $\ell_*(T)=\mathbb{E}\sup_{t\in T}|\sum_{i=1}^pg_it_i|$, where $t=(t_i)_{i=1}^p\in\mathbb{R}^p$.
\ede

Now we state following result which is a restricted eigenvalue condition for random matrices with independent entries where each entry is an isotropic $\psi_2$ probability measure:
\btheos[\cite{mendelson2007reconstruction}, Theorem 2.3]
\label{theorem:mendelson}
There exist absolute constants $c,\bar{c}>0$ for which the following holds. Let $T\subset\mathbb{R}^p$ be a star-shaped set and put $\mu$ to be an isotropic $\psi_2$ probability measure with constant $\alpha\geq1$. For $n\geq1$ let $X_1,...,X_n$ be independent, distributed according to $\mu$ and define $\Gamma=\sum_{i=1}^n\langle X_i,~.~\rangle e_i$, where $e_i\in\mathbb{R}^n$ is a vector with $i^{\text{th}}$-location to be 1 and all the other locations to be 0. If $0<f<1$, then with probability at least $1-\exp(-\bar{c}f^2n/\alpha^4)$, for all $x\in T$ such that $\|x\|_2\geq r_n^*(f/c\alpha^2)$, we have
\begin{equation}
\label{equation:RE}
(1-f)\|x\|_2\leq\frac{\|\Gamma x\|_2}{\sqrt{n}}\leq(1+f)\|x\|_2,
\end{equation}
where
$$r_n^*(f)=r_n^*(f,T):=\inf\{\rho>0: \rho\geq\ell_*(T_\rho)/(f\sqrt{n}) \}$$
and
$$T_{\rho}=\{x\in T;~\|x\|_2\leq\rho\}.$$
\etheos
Next we want to prove the restricted eigenvalue condition for subgaussian random matrices by using this theorem.
\spro
We first note that it is sufficient to prove this theorem for $\|x\|_2=1$. In fact if $x=\textbf{0}\in\mathbb{R}^p$ Theorem \ref{Theorem:RE} holds trivially. Otherwise we can consider the re-scaled vector $\tilde{x}=x/\|x\|_2$ with $\|\tilde{x}\|_2=1$. It can be seen that if this theorem holds for the re-scaled vector $\tilde{x}$, it also holds for $x$.

Next we define the set $V(r):=\{x\in\mathbb{R}^p~|~\|x\|_2=1,~\|x\|_1\leq r\}$, for a fixed radius $r>0$. It is possible that this set is empty for some choices of $r>0$, but we are only concerned with those choices for which it is non-empty. Define the random variable:
$$M(r,\Gamma):=\sup_{x\in V(r)}\{1-\frac{\|\Gamma x\|_2}{\sqrt{n}}\}.$$

Our goal is to show that with probability no larger than $\exp(-c_{\alpha}nf(r)^2)$ that 
$$M(r,\Gamma)=\sup_{x\in V(r)}\{1-\frac{\|\Gamma x\|_2}{\sqrt{n}}\}\geq\frac{3f(r)}{2},$$
where $f(r)=\frac{1}{4}+3c\alpha^2r\sqrt{\frac{\log p}{n}}$, $c_{\alpha}$ and $c$ are positive constants.

To see this,  we choose $f=f(r)=\frac{1}{4}+3c\alpha^2r\sqrt{\frac{\log p}{n}}$ in Theorem \ref{theorem:mendelson} (this $c$ here is defined in Theorem \ref{theorem:mendelson}). We first bound $\ell_*(V(r))$ as follows:
$$\ell_*(V(r))\leq r\mathbb{E}\max_{1\leq i\leq p}|g_i|\leq 3r\sqrt{\log p}$$
by using known results on Gaussian maxima (\cite{LedTal91}, Equation (3.13)). Then for all $x\in V(r)$ we have 
$$\|x\|_2=1\geq\frac{3c\alpha^2r\sqrt{\frac{\log p}{n}}}{f(r)}\geq\frac{c\alpha^2\ell_*(V(r))}{f(r)\sqrt{n}}=\frac{\ell_*(V(r))}{\frac{f(r)}{c\alpha^2}\sqrt{n}},$$
and by Theorem \ref{theorem:mendelson} with probability at least $1-\exp(-\bar{c}f(r)^2n/\alpha^4)$ we have for all $x\in V(r)$
$$1-f(r)\leq\frac{\|\Gamma x\|_2}{\sqrt{n}}.$$
Hence with probability no larger than $\exp(-9\bar{c}nf(r)^2/4\alpha^4)$,
$$M(r,\Gamma)=\sup_{x\in V(r)}\{1-\frac{\|\Gamma x\|_2}{\sqrt{n}}\}\geq\frac{3f(r)}{2}.$$

The remainder of the proof will mainly follow steps in \cite{RasWaiYu10b} where we use a peeling technique to extend our result to hold for $x$'s that have an arbitrary radius. We define the event
$$\Upsilon:=\{\exists~x\in\mathbb{R}^p~\text{s.t}~\|x\|_2=1~\text{and}~(1-\|\Gamma x\|_2)/\sqrt{n}\geq 3f(\|x\|_1)\}.$$
To prove the main theorem, the next step is to show that there are positive constants $c',c''$ such that $P[\Upsilon]\leq c'\exp(-c''n)$. Now we follow the standard peeling technique (\cite{vandeGeer, Alex85}) and we state the following lemma which is stated and proven in \cite{RasWaiYu10b}.
\blems[\cite{RasWaiYu10b}, Lemma 3]
\label{lemma:peeling}
Suppose that $d(v;X)$ is a random objective function with $v\in\mathbb{R}^p$ and $X$ is some random vector. $h:\mathbb{R}^p\to\mathbb{R}_+$ is some non-negative and increasing constraint function. $g:\mathbb{R}\to\mathbb{R}_+$ is a non-negative and strictly increasing function. A is a non-empty set. Moreover we suppose $g(r)\geq u$ for all $r\geq0$, and there exists some constant $c>0$ such that for all $r>0$, we have the tail bound
$$\mathbb{P}(\sup_{v\in A,h(v)\leq r}d(v;X)\geq g(r))\leq2\exp(-ca_ng^2(r)),$$
for some $a_n>0$. Then we have
$$\mathbb{P}[\mathcal{E}]\leq\frac{2\exp(-4ca_nu^2)}{1-\exp(-4ca_nu^2)},$$
where
$$\mathcal{E}:=\{\exists v\in A~\mbox{such that}~d(v;X)\geq2g(h(v))\}.$$
\elems
In order to use this lemma we choose the sequence $a_n=n$ and the set $A=\{x\in\mathbb{R}^p~|~\|x\|_2=1\}$, moreover we set
$$d(x,\Gamma)=1-\|\Gamma x\|_2/\sqrt{n},~~h(x)=\|x\|_1,~~\text{and}~~g(r)=3f(r)/2.$$
Since $g(r)=\frac{3f(r)}{2}=\frac{3}{8}+\frac{9}{2}c\alpha^2r\sqrt{\frac{\log p}{n}}\geq\frac{3}{8}$ for all $r>0$ and is strictly increasing, so that Lemma \ref{lemma:peeling} is applicable with $u=\frac{3}{8}$. Thus by using Lemma \ref{lemma:peeling} we have that $P[\Upsilon^c]\geq 1-c'\exp(-c''n)$ for some numerical constants $c'$ and $c''$.

Then for all $x\in\mathbb{R}^p$ with $\|x\|_2=1$, conditioned on the event $\Upsilon^c$ we have
$$1-\frac{\|\Gamma x\|_2}{\sqrt{n}}\leq 3f(\|x\|_1)=\frac{3}{4}+9c\alpha^2\|x\|_1\sqrt{\frac{\log p}{n}},$$
then
$$\frac{\|\Gamma x\|_2}{\sqrt{n}}\geq\frac{1}{4}-9c\alpha^2\|x\|_1\sqrt{\frac{\log p}{n}},$$
which completes the proof.
\fpro

\section{Appendix}

\subsection{Proof of Lemma \ref{lemma:ADsubgaussian}}

The proof for this lemma is mainly based on results in \cite{vershynin2010introduction}. First by using Lemma 5.5 in \cite{vershynin2010introduction} we know that the definition of the $\psi_2$ norm in \cite{vershynin2010introduction} is equivalent with our definition up to some absolute constant. Then we use a similar technique for the proof of Lemma 5.24 in \cite{vershynin2010introduction}. For every $x=(x_1,...,x_{p-1})\in S^{p-2}$ we have
$$\|\langle X\bar{D},x \rangle\|_{\psi_2}^2=\|\langle X, x\bar{D}^{\top} \rangle\|_{\psi_2}^2\leq C\sum_{i=1}^p(x\bar{D}^{\top})_i^2\|X_i\|_{\psi_2}^2\leq C\max_{1\leq i\leq p}\|X_i\|_{\psi_2}^2,$$
where the first inequality comes from Lemma 5.9 in \cite{vershynin2010introduction} and we also used $\sum_{i=1}^p(x\bar{D}^{\top})_i^2=(x\bar{D}^{\top})(x\bar{D}^{\top})^{\top}=xx^{\top}=1$ since $x\in S^{p-2}$. Since $X=(X_1,...,X_p)$ is distributed according to $\mu$, we have shown that $\|\langle X\bar{D},x \rangle\|_{\psi_2}$ is bounded by some absolute constant for every $x\in S^{p-2}$. It is also easy to see that $X\bar{D}$ is isotropic since $X$ is isotropic and $\bar{D}$ is orthonormal.

\subsection{Proof of Lemma \ref{lemma:lambdan}}

To bound $\|\frac{2n}{T}(y-TAD\theta^*)^{\top}A\bar{D}\|_{\infty}$ we need the following two lemmas. Lemma \ref{lemma:PoissonMGF} gives an upper bound for Poisson moment generating function and Lemma \ref{lemma:ADcondition} gives an upper bound for $\sum_{i=1}^n(AD)_{ij}^2$, with $2\leq j\leq p$.

\blems
\label{lemma:PoissonMGF}
If $W\sim\text{Poisson}(\lambda)$, for any $s\in [0,1]$ we have
\begin{eqnarray*}
\max\{\mathbb{E}\exp(s(W-\lambda)),\mathbb{E}\exp(s(\lambda-W))\}\leq\exp(\lambda s^2).
\end{eqnarray*}
\elems
\spro
By the moment generating function for Poisson random variable we know that
\begin{eqnarray*}
\mathbb{E}\exp(sW)=\exp(\lambda(e^s-1)),
\end{eqnarray*}
then we have
\begin{equation*}
\mathbb{E}\exp(s(W-\lambda))=\exp(\lambda(e^s-s-1)).
\end{equation*}
When $s\in [0,1]$ we have $e^s-s-1\leq s^2$, then we have shown
\begin{eqnarray*}
\mathbb{E}\exp(s(W-\lambda))\leq\exp(\lambda s^2).
\end{eqnarray*}
On the other hand we know that 
\begin{equation*}
\mathbb{E}\exp(s(\lambda-W))=\exp(\lambda(e^{-s}+s-1)).
\end{equation*}
When $s\in [0,1]$ we have $e^{-s}+s-1\leq s^2$, then we have shown
\begin{eqnarray*}
\mathbb{E}\exp(s(\lambda-W))\leq\exp(\lambda s^2).
\end{eqnarray*}
This completes the proof of Lemma \ref{lemma:PoissonMGF}.
\fpro

\blems
\label{lemma:ADcondition}
Under Assumption \ref{as:rip} we have for $2\leq j\leq p$,
$$\sum_{i=1}^n(AD)_{ij}^2\leq\frac{(1+\delta_{\tilde{K}})}{4n(a_u-a_{\ell})^2}.$$
\elems
\spro
By the definition of $\bar{D}$ and the construction of $A$
$$A\bar{D}=(\frac{\tilde{A}+\frac{a_u-2a_{\ell}}{\sqrt{n}}\mathbbm{1}_{n\times p}}{2(a_u-a_{\ell})\sqrt{n}})\bar{D}=\frac{\tilde{A}\bar{D}}{2(a_u-a_{\ell})\sqrt{n}}.$$
Then for $2\leq j\leq p$ we choose $u=e_j\in\mathbb{R}^p$ with $j$-th location to be 1 and all others to be 0, by Assumption \ref{as:rip},
$$\sum_{i=1}^n(AD)_{ij}^2=\|A\bar{D}e_j^{'}\|^2_2=\frac{\|\tilde{A}\bar{D}e_j^{'}\|_2^2}{4(a_u-a_{\ell})^2n}=\frac{\|\tilde{A}De_j\|_2^2}{4(a_u-a_{\ell})^2n}\leq\frac{(1+\delta_{\tilde{K}})}{4n(a_u-a_{\ell})^2},$$
where $e_j'\in\mathbb{R}^{p-1}$ is just $e_j$ without the first element 0.
\fpro

Then we can return to the proof of Lemma~\ref{lemma:lambdan}, note that
\begin{eqnarray*}
P(\|\frac{2n}{T}(y-TAD\theta^*)^{\top}A\bar{D}\|_{\infty}>t) &=& P(\max_{2\leq j\leq p}|\frac{2n}{T}\sum_{i=1}^n(y_i-T(AD\theta^*)_i)(AD)_{ij}|>t)\\
&=&P(\max_{2\leq j\leq p}|\sum_{i=1}^n(y_i-T(AD\theta^*)_i)\frac{n(AD)_{ij}}{T}|>\frac{t}{2}).
\end{eqnarray*}
Next we will first bound 
\begin{eqnarray}
\label{eq:t2left}
P(|\sum_{i=1}^n(y_i-T(AD\theta^*)_i)\frac{n(AD)_{ij}}{T}|>\frac{t}{2})&=&P(\sum_{i=1}^n(y_i-T(AD\theta^*)_i)\frac{n(AD)_{ij}}{T}>\frac{t}{2})\\
\label{eq:t2right}
&+&P(\sum_{i=1}^n(y_i-T(AD\theta^*)_i)\frac{n(AD)_{ij}}{T}<-\frac{t}{2}),
\end{eqnarray}
for each $2\leq i\leq p$ and then using union bound to get the final upper bound.\\

To bound (\ref{eq:t2left}) we have 
\begin{eqnarray*}
P(\sum_{i=1}^n(y_i-T(AD\theta^*)_i)\frac{n(AD)_{ij}}{T}>\frac{t}{2})=P(\sum_{i=1}^n(y_i-T(AD\theta^*)_i)\sqrt{\frac{n(AD)^2_{ij}}{M}}>\frac{tT}{2\sqrt{nM}}).
\end{eqnarray*}
From Lemma \ref{lemma:ADcondition} we can see that $\frac{n(AD)^2_{ij}}{M}\leq1$, then for $0\leq s\leq 1$, by using Chernoff's bound and Lemma \ref{lemma:PoissonMGF} we have 
\begin{eqnarray*}
&&P(\sum_{i=1}^n(y_i-T(AD\theta^*)_i)\frac{n(AD)_{ij}}{T}>\frac{t}{2})\\
&\leq&\exp(-\frac{stT}{2\sqrt{nM}})\prod_{i=1}^n\mathbb{E}\exp(s(y_i-T(AD\theta^*)_i)\sqrt{\frac{n(AD)^2_{ij}}{M}})\\
&\leq&\exp(-\frac{stT}{2\sqrt{nM}})\prod_{i=1}^n\exp(\frac{s^2n(AD)^2_{ij}}{M}*T(AD\theta^*)_i).
\end{eqnarray*}
By Lemma \ref{lemma:Afbound} we know that $T(AD\theta^*)_i=T(Af^*)_i\leq\frac{T}{n}$ for all $1\leq i\leq n$, then
\begin{eqnarray*}
P(\sum_{i=1}^n(y_i-T(AD\theta^*)_i)\frac{n(AD)_{ij}}{T}>\frac{t}{2})&\leq&\exp(-\frac{stT}{2\sqrt{nM}})\prod_{i=1}^n\exp(\frac{s^2n(AD)^2_{ij}}{M}*\frac{T}{n})\\
&=&\exp(-\frac{stT}{2\sqrt{nM}}+\frac{s^2T}{n}\sum_{i=1}^n\frac{n(AD)^2_{ij}}{M})\\
(\text{by Lemma \ref{lemma:ADcondition}})~&\leq&\exp(-\frac{stT}{2\sqrt{nM}}+\frac{s^2T}{n}).
\end{eqnarray*}
If we set $s=\frac{t}{4}\sqrt{\frac{n}{M}}$, then
\begin{eqnarray*}
P(\sum_{i=1}^n(y_i-T(AD\theta^*)_i)\frac{n(AD)_{ij}}{T}>\frac{t}{2})\leq\exp(-\frac{t^2T}{8M}+\frac{t^2T}{16M})=\exp(-\frac{t^2T}{16M}).
\end{eqnarray*}
On the other hand for (\ref{eq:t2right}), similarly for $0\leq s\leq 1$ we have
\begin{eqnarray*}
&&P(\sum_{i=1}^n(y_i-T(AD\theta^*)_i)\frac{n(AD)_{ij}}{T}<-\frac{t}{2})\\
&=&P(\sum_{i=1}^n(T(AD\theta^*)_i-y_i)\sqrt{\frac{n(AD)^2_{ij}}{M}}>\frac{tT}{2\sqrt{nM}})\\
&\leq&\exp(-\frac{stT}{2\sqrt{nM}})\prod_{i=1}^n\mathbb{E}\exp(s(T(AD\theta^*)_i-y_i)\sqrt{\frac{n(AD)^2_{ij}}{M}})\\
(\text{by Lemma \ref{lemma:PoissonMGF}}) &\leq&\exp(-\frac{stT}{2\sqrt{nM}})\prod_{i=1}^n\exp(\frac{s^2n(AD)^2_{ij}}{M}*T(AD\theta^*)_i)\\
(\text{by Lemma \ref{lemma:ADcondition}}) &\leq&\exp(-\frac{stT}{2\sqrt{nM}}+\frac{s^2T}{n}).
\end{eqnarray*}
For the same choice for $s=\frac{t}{4}\sqrt{\frac{n}{M}}$ we have 
\begin{eqnarray*}
P(\sum_{i=1}^n(y_i-T(AD\theta^*)_i)\frac{n(AD)_{ij}}{T}<-\frac{t}{2})\leq\exp(-\frac{t^2T}{16M}).
\end{eqnarray*}
Then by combining two result together we have 
\begin{eqnarray*}
P(|\sum_{i=1}^n(y_i-T(AD\theta^*)_i)\frac{n(AD)_{ij}}{T}|>\frac{t}{2})\leq2\exp(-\frac{t^2T}{16M}).
\end{eqnarray*}
Then by the union bound we have
\begin{eqnarray*}
P(\max_{2\leq j\leq p}|\sum_{i=1}^n(y_i-T(AD\theta^*)_i)\frac{n(AD)_{ij}}{T}|>\frac{t}{2})\leq2\exp(-\frac{t^2T}{16M}+\log(p-1)).
\end{eqnarray*}
If we set $t=\sqrt{\frac{32M\log p}{T}}$ then we have 
\begin{eqnarray*}
P(\|\frac{2n}{T}(y-TAD\theta^*)^{\top}A\bar{D}\|_{\infty}>\sqrt{\frac{32M\log p}{T}})\leq\frac{2}{p-1}.
\end{eqnarray*}
Also for this choice of $t$ we have $s=\sqrt{\frac{2n\log p}{T}}\in[0,1]$ by assumption that $T>2n\log p$. Thus with probability at least $1-\frac{2}{p-1}$,
$$\|\frac{2n}{T}(y-TAD\theta^*)^{\top}A\bar{D}\|_{\infty}\leq\sqrt{\frac{32M\log p}{T}}.$$

This completes the proof for Lemma \ref{lemma:lambdan}.

\subsection{Proof of Lemma \ref{lemma:rsc}}

Since $\hat{\Delta}\in\{\Delta\in\mathbb{R}^{p-1}~|~\|\Delta_{S^c_\eta}\|_1\leq3\|\Delta_{S_\eta}\|_1+4\|\bar{\theta}^*_{S^c_\eta}\|_1\}$, we have
\begin{eqnarray*}
\|\hat{\Delta}\|_1 &\leq& 4\|\hat{\Delta}_{{S_{\eta}}}\|_1+4\|\bar{\theta}^*_{S^c_{\eta}}\|_1\\
&\leq& 4\sqrt{|S_{\eta}|}\|\hat{\Delta}\|_2+4\|\bar{\theta}^*_{S^c_{\eta}}\|_1\\
&\leq& 4\sqrt{R_q}\eta^{-q/2}\|\hat{\Delta}\|_2+4\|\bar{\theta}^*_{S^c_{\eta}}\|_1,
\end{eqnarray*}
where the third inequality follows from (\ref{equation:Seta}). Then by Assumption \ref{Assumption:L1RE},
\begin{eqnarray*}
\sqrt{n}\|A\bar{D}\hat{\Delta}\|_2 &\geq& k_1\|\hat{\Delta}\|_2-4k_2\sqrt{\frac{\log p}{n}}(\sqrt{R_q}\eta^{-q/2}\|\hat{\Delta}\|_2+\|\bar{\theta}^*_{S^c_{\eta}}\|_1)\\
&=& \|\hat{\Delta}\|_2(k_1-4k_2\sqrt{\frac{R_q\log p}{n}}\eta^{-q/2})-4k_2\sqrt{\frac{\log p}{n}}\|\bar{\theta}^*_{S^c_{\eta}}\|_1.
\end{eqnarray*}
By choosing $\eta=\sqrt{\frac{\log p}{T}}$, we have
\begin{eqnarray*}
\sqrt{n}\|A\bar{D}\hat{\Delta}\|_2\geq(k_1-4k_2\sqrt{\frac{R_q\log p}{n}}(\frac{\log p}{T})^{-q/4})\|\hat{\Delta}\|_2-4k_2\sqrt{\frac{\log p}{n}}\|\bar{\theta}^*_{S^c_{\eta}}\|_1.
\end{eqnarray*}
Recall that by Assumption \ref{Assumption:ntildek} and the definition of $\tilde{K}$ we have $n\geq c_0R_q(\frac{\log p}{T})^{-q/2}\log p$ with constant $c_0>0$, then the above inequality becomes
$$\sqrt{n}\|A\bar{D}\hat{\Delta}\|_2\geq c'k_1\|\hat{\Delta}\|_2-c''k_2\sqrt{\frac{\log p}{n}}\|\bar{\theta}^*_{S^c_{\eta}}\|_1,$$
where $c',c''>0$ are some constants. We can then find constants $c_1,c_2>0$ such that 
$$n\|A\bar{D}\hat{\Delta}\|^2_2\geq c_1k_1^2\|\hat{\Delta}\|^2_2-c_2k_2^2\frac{\log p}{n}\|\bar{\theta}^*_{S^c_{\eta}}\|_1^2.$$

\bibliographystyle{plainnat}

\bibliography{Biblio_LL,Biblio_GroupLasso}

\end{document}

%% file: AMS_commands
\theoremstyle{plain}

\newtheorem{theo}{Theorem}[section]

\newtheorem{lem}{Lemma}[section]
\newtheorem{prop}{Proposition}[section]
\newtheorem{cor}{Corollary}[section]

\theoremstyle{definition} 

\newtheorem{nota}{Notation}[section]
\newtheorem{de}{Definition}[section]
\newtheorem{exa}{Example}[section]
\newtheorem{as}{Assumption}[section]
\newtheorem{alg}{Algorithm}[section]
\newtheorem{rem}{Remark}[section]

\newcommand{\btheo}{\begin{theo}}
\newcommand{\bde}{\begin{de}}
\newcommand{\ble}{\begin{lem}}
\newcommand{\bpr}{\begin{prop}}
\newcommand{\bno}{\begin{nota}}
\newcommand{\bex}{\begin{exa}}
\newcommand{\bcor}{\begin{cor}}
\newcommand{\spro}{\begin{proof}}
\newcommand{\bas}{\begin{as}}
\newcommand{\balg}{\begin{alg}}
\newcommand{\brem}{\begin{rem}}

\newcommand{\etheo}{\end{theo}}
\newcommand{\ede}{\end{de}}
\newcommand{\ele}{\end{lem}}
\newcommand{\epr}{\end{prop}}
\newcommand{\eno}{\end{nota}}
\newcommand{\eex}{\end{exa}}
\newcommand{\ecor}{\end{cor}}
\newcommand{\fpro}{\end{proof}}
\newcommand{\eas}{\end{as}}
\newcommand{\ealg}{\end{alg}}
\newcommand{\erem}{\end{rem}}

\theoremstyle{plain}

\newtheorem{theos}{Theorem}
\newtheorem{props}{Proposition}
\newtheorem{lems}{Lemma}
\newtheorem{cors}{Corollary}

\theoremstyle{definition}
\newtheorem{exas}{Example}
\newtheorem{algs}{Algorithm}
\newtheorem{asss}{Asumption}
\newtheorem{defns}{Definition}

\newcommand{\btheos}{\begin{theos}}
\newcommand{\etheos}{\end{theos}}
\newcommand{\bprops}{\begin{props}}
\newcommand{\eprops}{\end{props}}
\newcommand{\bdes}{\begin{defns}}
\newcommand{\edes}{\end{defns}}
\newcommand{\blems}{\begin{lems}}
\newcommand{\elems}{\end{lems}}
\newcommand{\bcors}{\begin{cors}}
\newcommand{\ecors}{\end{cors}}
\newcommand{\bexs}{\begin{exas}}
\newcommand{\eexs}{\end{exas}}
\newcommand{\balgs}{\begin{algs}}
\newcommand{\ealgs}{\end{algs}}
\newcommand{\bass}{\begin{asss}}
\newcommand{\eass}{\end{asss}}

%% file: Poisson_Revision_v1.bbl
\begin{thebibliography}{10}

\bibitem{riceCamera}
M.~F. Duarte, M.~A. Davenport, D.~Takhar, J.~N. Laska, T.~Sun, K.~F. Kelly, and
  R.~G. Baraniuk.
\newblock Single pixel imaging via compressive sampling.
\newblock {\em IEEE Sig. Proc. Mag.}, 25(2):83--91, 2008.

\bibitem{CandesSIM}
V.~Studer, J.~Bobin, M.~Chahid, H.~S. Mousavi, E.~Candes, and M.~Dahan.
\newblock Compressive fluorescence microscopy for biological and hyperspectral
  imaging.
\newblock {\em Proceedings of the National Academy of Sciences of the United
  States of America (PNAS)}, 109(26), 2012.

\bibitem{nimmi_spie}
Z.~T. Harmany, J.~Mueller, Q.~Brown, N.~Ramanujam, and R.~Willett.
\newblock Tissue quantification in photon-limited microendoscopy.
\newblock In {\em Proc. SPIE Optics and Photonics}, 2011.

\bibitem{mca}
J.~Bobin, J.-L. Starck, J.~Fadili, Y.~Moudden, and D.~L. Donoho.
\newblock Morphological component analysis: An adaptive thresholding strategy.
\newblock {\em IEEE Transactions on Image Processing}, 16(11):2675--2681, 2007.

\bibitem{gmca}
J.~Bobin, J.-L. Starck, J.~Fadili, and Y.~Moudden.
\newblock Sparsity and morphological diversity in blind source separation.
\newblock {\em IEEE Transactions on Signal Processing}, 16(11):2662--2674,
  2007.

\bibitem{ElephantsMice}
C.~Estan and G.~Varghese.
\newblock New directions in traffic measurement and accounting: focusing on the
  elephants, ignoring the mice.
\newblock {\em ACM Trans. Computer Sys.}, 21(3):270--313, 2003.

\bibitem{CounterBraids}
Y.~Lu, A.~Montanari, B.~Prabhakar, S.~Dharmapurikar, and A.~Kabbani.
\newblock Counter {B}raids: a novel counter architecture for per-flow
  measurement.
\newblock In {\em Proc. ACM SIGMETRICS}, 2008.

\bibitem{fishing}
M.~Raginsky, S.~Jafarpour, R.~Willett, and R.~Calderbank.
\newblock Fishing in {P}oisson streams: focusing on the whales, ignoring the
  minnows.
\newblock In {\em Proc. Forty-Fourth Conference on Information Sciences and
  Systems}, 2010.
\newblock \href{http://arxiv.org/abs/1003:2836}{arXiv:1003.2836}.

\bibitem{expander_pcs}
M.~Raginsky, S.~Jafarpour, Z.~Harmany, R.~Marcia, R.~Willett, and
  R.~Calderbank.
\newblock Performance bounds for expander-based compressed sensing in {P}oisson
  noise.
\newblock {\em IEEE Transactions on Signal Processing}, 59(9), 2011.
\newblock \href{http://arxiv.org/abs/1007.2377}{arXiv:1007.2377}.

\bibitem{LaureDNA}
L.~Sansonnet.
\newblock Wavelet thresholding estimation in a {Poisson}ian interactions model
  with application to genomic data.
\newblock {\em Scandinavian Journal of Statistics}, 2013.
\newblock doi: 10.1111/sjos.12009.

\bibitem{dense_pcs}
M.~Raginsky, R.~Willett, Z.~Harmany, and R.~Marcia.
\newblock Compressed sensing performance bounds under {P}oisson noise.
\newblock {\em IEEE Transactions on Signal Processing}, 58(8):3990--4002, 2010.
\newblock \href{http://arxiv.org/abs/0910.5146}{arXiv:0910.5146}.

\bibitem{jiang2014minimax}
X.~Jiang, G.~Raskutti, and R.~Willett.
\newblock Minimax optimal rates for poisson inverse problems with physical
  constraints.
\newblock {\em arXiv preprint arXiv:1403.6532}, 2014.

\bibitem{willett2011poisson}
R.~M. Willett and M.~Raginsky.
\newblock {Poisson} compressed sensing.
\newblock {\em Defense Applications of Signal Processing}, 2011.

\bibitem{foucart2010gelfand}
Simon Foucart, Alain Pajor, Holger Rauhut, and Tino Ullrich.
\newblock The {G}elfand widths of $\ell_p$-balls for $0<p\leq1$.
\newblock {\em Journal of Complexity}, 26(6):629--640, 2010.

\bibitem{RasWaiYu11}
G.~Raskutti, M.~J. Wainwright, and B.~Yu.
\newblock Minimax rates of estimation for high-dimensional linear regression
  over $\ell_q$-balls.
\newblock {\em IEEE Transactions of Information Theory}, 57(10):6976--6994,
  2011.

\bibitem{Neg10}
S.~Negahban, P.~Ravikumar, M.~J. Wainwright, and B.~Yu.
\newblock A unified framework for high-dimensional analysis of $m$-estimators
  with decomposable regularizers.
\newblock {\em Statistical Science}, 27(4):538--557, 2012.

\bibitem{bobkov1998modified}
S.~G. Bobkov and M.~Ledoux.
\newblock On modified logarithmic sobolev inequalities for bernoulli and
  poisson measures.
\newblock {\em Journal of Functional Analysis}, 156(2):347--365, 1998.

\bibitem{Baraniuk08}
R.~Baraniuk, M.~Davenport, R.~DeVore, and M.~Wakin.
\newblock A simple proof of the restricted isometry property for random
  matrices.
\newblock {\em Constructive Approximation}, 28(3):253--263, 2008.

\bibitem{Kuh01}
T.~K\"{u}hn.
\newblock A lower estimate for entropy numbers.
\newblock {\em Journal of Approximation Theory}, 110:120--124, 2001.

\bibitem{BicRitTsy08}
P.~Bickel, Y.~Ritov, and A.~Tsybakov.
\newblock Simultaneous analysis of {L}asso and {D}antzig selector.
\newblock Submitted to Annals of Statistics, 2008.

\bibitem{RasWaiYu10b}
G.~Raskutti, M.~J. Wainwright, and B.~Yu.
\newblock Restricted eigenvalue conditions for correlated {G}aussian designs.
\newblock {\em Journal of Machine Learning Research}, 11:2241--2259, 2010.

\bibitem{GeerBuhl09}
S.~van~de Geer and P.~Buhlmann.
\newblock On the conditions used to prove oracle results for the lasso.
\newblock {\em Electronic Journal of Statistics}, 3:1360--1392, 2009.

\bibitem{mendelson2007reconstruction}
S.~Mendelson, A.~Pajor, and N.~Tomczak-Jaegermann.
\newblock Reconstruction and subgaussian operators in asymptotic geometric
  analysis.
\newblock {\em Geometric and Functional Analysis}, 17(4):1248--1282, 2007.

\bibitem{zhou2009restricted}
S.~Zhou.
\newblock Restricted eigenvalue conditions on subgaussian random matrices.
\newblock {\em arXiv preprint arXiv:0912.4045}, 2009.

\bibitem{Jiang2015data}
X.~Jiang, P.~Reynaud-Bouret, V.~Rivoirard, L.~Sansonnet, and R.~Willett.
\newblock A data-dependent weighted lasso under poisson noise.
\newblock {\em arXiv preprint arXiv:1509.08892}, 2015.

\bibitem{harmany2012spiral}
Zachary~T Harmany, Roummel~F Marcia, and Rebecca~M Willett.
\newblock This is spiral-tap: Sparse poisson intensity reconstruction
  algorithms�theory and practice.
\newblock {\em IEEE Transactions on Image Processing}, 21(3):1084--1096, 2012.

\bibitem{han1994generalizing}
T.~S. Han and S.~Verdu.
\newblock Generalizing the fano inequality.
\newblock {\em IEEE Transactions on Information Theory}, 40(4):1247--1251,
  1994.

\bibitem{IbrHas81}
I.~A. Ibragimov and R.~Z. Has'minskii.
\newblock {\em Statistical {E}stimation: {A}symptotic {T}heory}.
\newblock Springer-Verlag, New York, 1981.

\bibitem{YanBar99}
Y.~Yang and A.~Barron.
\newblock Information-theoretic determination of minimax rates of convergence.
\newblock {\em Annals of Statistics}, 27(5):1564--1599, 1999.

\bibitem{BiRiTsy08}
P.~Bickel, Y.~Ritov, and A.~Tsybakov.
\newblock Simultaneous analysis of {L}asso and {D}antzig selector.
\newblock {\em Annals of Statistics}, 37(4):1705--1732, 2009.

\bibitem{vandeGeer}
S.~van~de Geer.
\newblock {\em Empirical Processes in M-Estimation}.
\newblock Cambridge University Press, 2000.

\bibitem{LedTal91}
M.~Ledoux and M.~Talagrand.
\newblock {\em Probability in Banach Spaces: Isoperimetry and Processes}.
\newblock Springer-Verlag, New York, NY, 1991.

\bibitem{Alex85}
K.~S. Alexander.
\newblock Rates of growth for weighted empirical processes.
\newblock In {\em Proceedings of the Berkeley Conference in Honor of Jerzy
  Neyman and Jack Kiefer}, pages 475--493, Berkeley, 1985. UC Press.

\bibitem{vershynin2010introduction}
R.~Vershynin.
\newblock Introduction to the non-asymptotic analysis of random matrices.
\newblock {\em arXiv preprint arXiv:1011.3027}, 2010.

\end{thebibliography}
